\documentclass{article}

\usepackage{amsmath, amsthm, amsfonts, amssymb}

\usepackage[utf8]{inputenc}
\usepackage[T1]{fontenc}
\usepackage{enumitem}

\date{}

\newtheorem{thm}{Theorem}
\newtheorem{cor}[thm]{Corollary}
\newtheorem{lem}[thm]{Lemma}
\newtheorem{prop}[thm]{Proposition}
\theoremstyle{definition}
\newtheorem{defn}[thm]{Definition}
\theoremstyle{remark}

\def\R{\mathbb{R}}

\def\C{\mathbb{C}}
\def\B{\mathbb{B}}
\def\D{\mathbb{D}}

\def\exp{\, \mbox{exp}}

\def\dist{\, \mbox{\normalfont{dist}}}

\def\Jet{\, \mbox{\normalfont{Jet}}}
\def\Par{\,  \widetilde{ \mbox{Jet} }}
\def\Ric{\, \mbox{\normalfont{Ric}}}
\def\Scal{\, \mbox{\normalfont{Scal}}}
\def\proof{{\it Proof. }}
\def\codim{{\rm codim }}
\def\Gdec{G_{\rm dec}}
\def\Pdec{P_{\rm dec}}
\def\II{\mbox{II}}

\def\cL{\cal{L}}

\def\cY{{\cal Y}}

\def\Pol{{\cal P}ol}
\def\supp{{\rm supp}}

\def\vol{{\rm vol}}

\title{Construction of negatively curved complete intersections}
\author{Jean-Paul Mohsen\thanks{The author is supported by the UMI-CRM 3457 of the CNRS in Montreal.}}
\begin{document}
\maketitle

\begin{abstract}
Using the Donaldson-Auroux theory, 
we construct complete intersections 
in complex projective manifolds, which are negatively curved in various ways.
In particular, we prove the existence of compact simply connected K\"ahler manifolds with negative 
holomorphic bisectional curvature.
We also construct hyperbolic hypersurfaces and we obtain bounds for their Kobayashi hyperbolic metric.
\end{abstract}

In this article, we give new applications of the asymptotic methods invented by Donaldson in \cite{Do96}.
Although the theory was primarily designed to prove structure theorems in symplectic geometry,
applications to complex geometry already appear in 
\cite[Corollary 33]{Do96}.
Nonetheless, it seems that for a number of years
complex applications haven't been continued.
Recently, Giroux and Pardon have used Donaldson's techniques in their study of Stein manifolds
\cite{GiPa17}.
In this paper we give other applications in the context of complex projective geometry.

Let $X$ be a complex projective manifold endowed with an ample line bundle $L \rightarrow X$.
A hypersurface of degree $k$ in $X$ is, by definition, the zero set $Y = \{ s_k=0 \}$ of a holomorphic section $s_k$ of the line bundle $L^{\otimes k}$
(the $k^{\mbox{\footnotesize{th}}}$ tensor power of $L$). More generally,
a complete intersection $Y\subset X$ of dimension $d$ defined by equations of degree $k$ is the zero set $ Y = \{ s_k=0 \}$ of a holomorphic section $s_k$ 
of $\C^{n-d} \otimes L^{\otimes k}=L^{\otimes k} \oplus \dots \oplus L^{\otimes k}$ where $n = \dim X$.
In this paper, we will always assume that $Y$ is smooth (with $s_k$ transverse to $0$).

Fix a Hermitian metric $\mu$ on $X$. Our main results are existence results for complete intersections 
in $X$ whose curvature for the induced metric has various negativity properties.


\begin{thm}\label{negatively curved}

Let $X$ be a complex projective manifold of dimension $n$
equipped with a Hermitian metric $\mu$ and an ample line bundle $L \rightarrow X$.
\begin{enumerate}[label=(\alph*)]
\item \label{curve}
If $n\geq 3$
then, for every sufficiently large $k$, there exists a curve $Y \subset X$ 
which is a complete intersection  defined by equations of degree $k$ 
with negative curvature (for the induced metric).

\item\label{item Ricci} 
If $d\leq n - 2$
then, for every sufficiently large $k$, there exists a complete intersection $Y$ 
of dimension $d$  defined by equations of  degree $k$
with negative Ricci curvature in $X$.

\item\label{item Scal}  
If $d\leq n - 1$ and $n\geq 3$
then, for every sufficiently large $k$, there exists a complete intersection $Y$ 
of dimension $d$ defined by equations of  degree $k$
with negative scalar curvature in $X$.

\item\label{item holsec} 
If $n\geq 3d$
then, for every sufficiently large $k$, there exists a complete intersection $Y$ 
of dimension $d$ defined by equations of  degree $k$
with negative holomorphic sectional curvature in $X$.

\item\label{item bisec}
If $n\geq 4d - 1$
then, for every sufficiently large $k$, there exists a complete intersection $Y$ 
of dimension $d$ defined by equations of  degree $k$
with negative holomorphic bisectional curvature in $X$.

\end{enumerate}
\end{thm}

\bigskip

We could find no trace of such results in the literature, so we presume they are new.

\bigskip

Remark.
In fact, we prove a stronger result. In each case, we construct a 
sequence of submanifolds $(Y_k)_{k \gg 1}$ independent of the metric $\mu$, such that, for every $k$, the submanifold $Y_k \subset X$ is
a complete intersection of  dimension $d$  defined by equations of degree $k$
and, 
for every metric $\mu$, if $k$ is sufficiently large (in a sense which depends on $\mu$), then the corresponding curvature of $Y_k$ is negative.
In particular, if $\mu_1$, $\mu_2$, ... are a finite family of Hermitian metrics then there exists a complete 
intersection whose curvature is negative for every metric in the family.
The same holds for a compact family of Hermitian metrics.

\medskip

A simple example shows that 
Theorem \ref{negatively curved}\ref{curve} is sharp: 
no curve in the projective plane $\C P^2$ can have nonpositive curvature everywhere for the metric induced by the Fubini-Study metric. 
Similarly, if $X$ is an abelian variety with a constant K\"ahler metric then no
complex hypersurface $Y \subset X$
can have negative Ricci curvature everywhere and therefore Theorem \ref{negatively curved}\ref{item Ricci} is sharp.
For $n\geq 2d$, Brotbek, Darondeau and Xie (\cite{BrDa}, \cite{Xi}) have constructed complete intersections with 
ample cotangent bundles. 
In a way,
Theorem \ref{negatively curved}\ref{item bisec} is a metric version of this algebraic result, under the stronger hypothesis $n\geq 4d-1$.  

Theorem \ref{negatively curved}\ref{item bisec} answers a classical question concerning the bisectional curvature (see \cite{Wo81} and \cite[Question 35]{Ya82}).
If $X$ is simply connected then the Lefschetz theorem implies that complete intersections of dimension $d \geq 2$ are simply connected.
Hence, Theorem \ref{negatively curved}\ref{item bisec} implies the following result:

\begin{cor}\label{negative bisectional}
There exist compact simply connected K\"ahler manifolds with negative 
holomorphic bisectional curvature.
\end{cor}
In particular, this corollary shows that a simply connected complete K\"ahler manifold 
with negative holomorphic bisectional curvature need not be Stein, which was not known
(\cite[Question 35]{Ya82}).

The existence of compact K\"ahler surfaces 
with negative holomorphic bisectional curvature and $b_1=0$ 
is well known to experts (e.g. Mumford's fake projective plane \cite{Mu79}).
\bigskip

If $Y\subset X$ is a complex submanifold then the metric $\mu$ induces Hermitian metrics on the (complex) cotangent bundle $T^{*} Y$,
on the normal bundle $NY$ and, more generally, on their exterior powers $\bigwedge^l (T^{*} Y)$ and $\bigwedge^l (N Y)$.
Our next result constructs submanifolds for which the Hermitian bundle $\bigwedge^l (T^{*} Y)$ (or $\bigwedge^l (N Y)$) is Griffiths-positive.

\begin{thm}\label{exterior}
Let $X$ be a complex projective manifold of dimension $n$
equipped with a Hermitian metric $\mu$ and an ample line bundle $L \rightarrow X$.

\begin{enumerate}[label=(\alph*)]
\item\label{puissances cotangent}
If $2d (1+l) \leq l (n+l)$
then, for every sufficiently large $k$, there exists a complete intersection $Y$ of dimension $d$  defined by equations of  degree $k$ 
such that the Hermitian bundle $\bigwedge^l (T^{*} Y)$ is Griffiths-positive (for the induced metric).

\item\label{puissances normal} 
Similarly,
if $l (n - l) \leq 2d ( l - 1 )$ then there exists $Y$
such that $\bigwedge^l (N Y)$ is Griffiths-positive.
\end{enumerate}
\end{thm}
Remark. 
Theorem \ref{exterior}\ref{puissances cotangent} with $l=1$ is just Theorem \ref{negatively curved}\ref{item bisec}, if the metric $\mu$ is K\"ahler. 
Indeed, the cotangent bundle is Griffiths-positive if and only if the tangent bundle is Griffiths-negative, in other words if and only if the holomorphic bisectional curvature is negative
(see definitions in Sections \ref{proof bisec} and \ref{def Griffiths} below).
Clearly, if $l=1$ then Theorem \ref{exterior}\ref{puissances normal} is empty.
Furthermore, Theorem \ref{exterior}\ref{puissances cotangent} with $l=d$ is a reformulation of Theorem \ref{negatively curved}\ref{item Ricci}, if the metric $\mu$ is K\"ahler. 

\bigskip

Masuda and Noguchi \cite{MaNo96} have proved that every projective manifold contains hyperbolic hypersurfaces
(in the sense of Brody or Kobayashi).
Theorem \ref{hyperbo} below
gives a new construction of such hyperbolic hypersurfaces whose Kobayashi hyperbolic metric actually satisfies quantitative bounds. 
(Recently, it has been shown that almost every 
hypersurface of sufficiently large degree
is hyperbolic, see \cite{Br17} and \cite{Si15}.)

\begin{thm}\label{hyperbo}
Let $X$ be a complex projective manifold
equipped with a Hermitian metric $\mu$ and an ample line bundle $L \rightarrow X$ 
and let $\D \subset \C$ denote the unit disk.

Then, for every sufficiently large $k$, there exists a hyperbolic hypersurface $Y$ of degree $k$ such that every holomorphic map $f: \D \rightarrow Y$
satisfies
$$
\| f'(0) \| \leq \frac{C}{\sqrt{k}}
$$
for some constant $C$ which only depends on $X$, $L$ and $\mu$.
\end{thm}

\medskip

The organization of the paper is as follows. 
The Donadson-Auroux theory is discussed in Section \ref{att}. Results are stated in the language of \cite{Mo}.
In Sections \ref{preuve type}, \ref{proof exterior} and \ref{proof hyperbo}, we use the theory of Donaldson and Auroux in order to prove Theorem \ref{negatively curved}, Theorem \ref{exterior}
and Theorem \ref{hyperbo}.
In Section \ref{la theorie}, we prove the main result of Section \ref{att}.

\medskip

{\bf Acknowledgements.}
Emmanuel Giroux's ideas are the main inspiration for the theory which I present here.
I thank him for many crucial ideas and it's a pleasure to say that
all topics in this paper have been discussed with him for decades. 
I would like to thank Lionel Darondeau and Steven Lu for useful discussions concerning applications of asymptotic techniques to complex geometry.
In particular, the idea of applying these techniques in order to get Corollary \ref{negative bisectional} is due to Steven Lu.
I thank Misha Kapovich for references.
I also thank the UMI-CRM 3457 of the CNRS in Montreal for 
their hospitality during the preparation of this manuscript.

\section{Asymptotic transversality theory}\label{att}
\subsection{Limit submanifolds}\label{limit sub}
Let $X$ be a complex projective manifold of dimension $n$ equipped with an ample line bundle $L \rightarrow X$.
Recall that a submanifold $Y$ is a complete intersection of dimension $d$ defined by equations of degree $k$
if it is the common zero set of $n - d$ holomorphic sections of $L^{\otimes k}$.

In \cite{Do96} and \cite{Au97}, Donaldson and Auroux construct complete intersections $Y_k$
of fixed dimension $d$  defined by equations of degree $k$ going to $\infty$ which satisfy some remarkable compactness properties.
More specifically, the geometry of $Y_k$, when regarded at scale $\frac{1}{\sqrt{k}}$, remains bounded
(see Definition \ref{renorm sequence} below),
in particular the submanifolds $(Y_k)_{k\gg 1}$ admit limits which are submanifolds in $\C^n$
just because a complex manifold viewed at smaller and smaller scale looks like $\C^n$.
We will give below some precise definitions and results in which we try to emphasize the role of this limit submanifolds.

\medskip

Let's give formal definitions.
For every sufficiently large $k$, let $Y_k \subset X$ be a smooth complete intersection of dimension $d$ defined by equations of degree $k$,
let $\B = \B(1) \subset \C^n$ be the unit ball and, more generally, let $\B(r)$ denote the ball of radius $r$.
Given holomorphic charts $\varphi_k: \B \rightarrow X$, for $k$ in some infinite set $I$ of integers, we define the corresponding 
{\it rescaled charts} $R\varphi_k: \B(\sqrt{k}) \rightarrow X$ by the following formula:
$$
R\varphi_k(p)  = \varphi_k \left( \frac{p}{\sqrt{k}}  \right)
$$
where the point $p$ belongs to $\B(\sqrt{k})$.
For every sufficently large $k \in I$, we define the {\it renormalized submanifold} $RY_k = (R\varphi_k)^{-1} (Y_k) \subset \C^n$.

This is a rather abbreviated notation. 
It is important to remember that the submanifold $RY_k$ depends on the chart $\varphi_k$, even though the chart doesn't appear in the notation.

\begin{defn}\label{renorm sequence}
The sequence $(Y_k)_{k\gg 1}$ is {\it renormalizable} if for every sequence $(\varphi_k: \B \rightarrow X)_{k\in I}$ of holomorphic charts (indexed 
by some infinite set $I$ of sufficiently large integers) which satisfies properties 1 and 2 below, 
the corresponding family $(RY_k)_{k\in I}$ of renormalized submanifolds is
relatively compact for the smooth compact-open topology. The charts $(\varphi_k: \B \rightarrow X)_{k\in I}$
are supposed to satisfy the following properties:

\begin{enumerate}
\item
Normality: The family $(\varphi_k)_{k\in I}$ is normal ({\it i.e.} relatively compact for the smooth compact-open topology on holomorphic maps). 

\item
Non degeneracy:
If $\varphi_{\infty}: \B \to X$ is a limit map of the normal sequence $(\varphi_k)_{k\in I}$
then the tangent map $d \varphi _{\infty} (0)$
is an isomorphism.

\end{enumerate}

Note that we don't assume that the points $(\varphi_k (0))_{k\in I}$ form a constant sequence in $X$.

Let $(RY_k)_{k\in J}$ be a subsequence 
which converges to some submanifold
$Y_{\infty} \subset \C^n$.
Then $Y_{\infty}$ is a complex submanifold of dimension $d$ (in the setting of complex submanifolds, the smooth convergence is the same 
as the holomorphic convergence).
Such a $Y_{\infty}$ is called a {\it limit submanifold} of the renormalizable sequence $(Y_k)_{k\gg 1}.$
\end{defn}

\medskip
The following theorem is a reformulation in the above terminology of results due to Donaldson  \cite{Do96}
and Auroux  \cite{Au97}.

\medskip
\begin{thm}\label{diviseur de Donaldson}
Let $X$ be a complex projective manifold of dimension $n$ equipped with an ample line bundle $L$ and
let $d$ be an integer with $1 \leq d < n$.

Then, for every sufficiently large $k$, there exists a complete intersection $Y_k \subset X$ of dimension $d$ defined by equations of degree $k$
such that submanifolds $(Y_k)_{k\gg 1}$ form a renormalizable sequence.
\end{thm}

This Theorem is a special case of Theorem \ref{avoidance} below.

\subsection{Jets of submanifolds}\label{jets}
Let $\Jet^l_{d,n}$ denote the space of $l-$jets of complex submanifolds in $\C^n$ of dimension $d$.
In order to describe $\Jet^l_{d,n}$, consider the set $\Par^l_{d,n}$ of  $l-$jets of regular parametizations of submanifolds.
Hence, the space $\Par^l_{d,n}$ is a (smooth, quasi-projective) algebraic manifold whose elements are $l-$jets of (germs of) holomorphic 
maps $(\C^d, 0) \to \C^n$
whose differential at $0$ is injective.
In particular, it is endowed with the Zariski topology.
Note that: $\Jet^l_{d,n} = \Par^l_{d,n} / G_d^l$ where $G_d^l$ is the group
of $l-$jets of (germs of) biholomorphisms $(\C^d, 0)  \to  (\C^d, 0)$.
In this paper $\Jet^l_{d,n}$ is endowed with the quotient topology and closed subsets in $\Jet^l_{d,n}$ are called closed {\it complex algebraic} subsets.
This terminology comes from the fact that 
one can endow $\Jet^l_{d,n}$ with a natural structure of algebraic manifold and the quotient topology above
is just the corresponding Zariski topology.  

\medskip

Notice $\Jet^{l+1}_{d,n} \rightarrow \Jet^l_{d,n}$ is a fiber bundle.
Let $Y\subset \C^n$ be a complex submanifold of dimension $d$ and let $p\in Y$ be a point.
Of course, 
$\Jet^0_{d,n}=\C^n$
since the $0-$jet of $Y$ at $p$ is just the point $p$. 
The $1-$jet of $Y$ at $p$ is a pair $(p,T)$ where $T\subset \C^n$ is a linear subspace (tangent space).
Denote by Grass($d,\C^n)$ the grassmannian of linear subspaces of dimension $d$ in $\C^n$.
Then $\Jet^1_{d,n} = \C^n \times$Grass$(d,\C^n)$.

The $2-$jet of $Y$ at $p$ is a triple $(p,T,\II)$ where $T$ is a linear subspace and $\II$ is a symmetric bilinear map
$T\times T \rightarrow \C^n / T$ (the second fundamental form).
Hence, $\Jet^2_{d,n} \rightarrow \Jet^1_{d,n}$ is a vector bundle having rank $\frac{d(d+1)(n-d)}{2}$.
The fiber over $(p,T)\in \Jet^1_{d,n}$ equals the space Sym$^2(T, \C^n / T)$ of symmetric bilinear maps from 
$T\times T$ to $\C^n / T$. 
If $a\in \Jet^2_{d,n}$ then we will use the notation $a = (p_a, T_a , \II _a)$.

\subsection{Avoidance theorem}
Consider a subset $A$ in $\Jet^l_{d,n}$ satisfying the following two conditions:
\begin{enumerate}
\item
$A$ is a closed complex algebraic subset

\item
$A$ is invariant under the natural action of affine transformations of $\C^n$ upon $\Jet^l_{d,n}$.

\end{enumerate}

Our main tool is a generalization of Theorem \ref{diviseur de Donaldson}:
\begin{thm}\label{avoidance}
Let $X$ be a complex projective manifold of dimension $n$ equipped with an ample line bundle $L$ and
let $A \subset \Jet^l_{d,n}$ be a subset satifying the above conditions 1 and 2.

If $\codim (A) > d$
then, for every sufficiently large $k$, there exists a complete intersection $Y_k \subset X$ of dimension $d$ defined by equations of degree $k$
such that submanifolds $(Y_k)_{k\gg 1}$ form a renormalizable sequence
which {\rm avoids the subset $A$ asymptotically}
in the sense that the $l-$jets of its limit submanifolds lie in the complement $\Jet^l_{d,n}\backslash A$.
\end{thm}

The proof is an application of Donaldson and Auroux's techniques. 
We will give a detailed proof of Theorem \ref{avoidance} in Section \ref{la theorie}.

\medskip

In order to apply the avoidance theorem, we need to produce large codimension subsets.
That's the use of the {\it main theorem of elimination theory.} 

\begin{thm}\label{proj codim}
Let $X_1$, $X_2$ be complex algebraic manifolds.
Assume $X_1$ is quasi-projective and $X_2$ is projective.
 Let $A \subset X_1 \times X_2$ be a closed algebraic subset
and denote by $\pi_1$ the first projection.

Then $\pi_1 (A)$ is a closed algebraic subset in $X_1$ with 
$\dim (\pi_1 (A)) \leq \dim ( A )$ 
(hence: $\codim  (\pi_1 (A)) \geq \codim ( A ) - \dim ( X_2 )$).
 
\end{thm}

This classical result is proved in \cite[Section I.9]{Mu99}.

\section{Proof of Theorem \ref{negatively curved}}\label{preuve type}

As we shall see, all cases in Theorem \ref{negatively curved} are proved in the same way.

\medskip

Let $X$ be a complex projective manifold of dimension $n$
equipped with a Hermitian metric $\mu$ and an ample line bundle $L \rightarrow X$.
All types of curvature in Theorem \ref{negatively curved} are non-positive for a complex submanifold $Y \subset \C^n$.
(Here, $\C^n$ is endowed with a constant K\"ahler metric and $Y$ is endowed with the induced metric.)
For each type of curvature, we will construct a closed algebraic subset 
$A \subset \Jet^2_{d,n}$ of codimension $>d$
which satisfies the invariance assumption of Theorem \ref{avoidance}
and the following assumption:
\medskip

{\it 
If $Y \subset \C^n$ is a complex submanifold of dimension $d$ and if $p\in Y$ is a point
such that the $2-$jet of $Y$ at $p$ belongs to the complement 
$\Jet^2_{d,n} \backslash A$
then the curvature of $Y$ at $p$ is negative.
}
\medskip

If such an $A$ exists then apply Theorem \ref{avoidance}.
This theorem provides a renormalizable sequence of complete intersections $(Y_k)_{k\gg 1}$ in $X$
of dimension $d$
which avoids the subset $A$ asymptotically.
Then the limit submanifolds $Y_{\infty}$ have negative curvature (for every constant K\"ahler metric on $\C^n$ since
the subset $A$ satisfies the invariance assumption).
Let's prove that, for every sufficiently large $k$, the submanifold $Y_k$ has negative curvature.

Our proof goes by contradiction. Assume that, for infinitely many integers $k$, there exists some point $p_k \in Y_k$
where the curvature of $Y_k$ {\it isn't} negative. This holds for both metrics $\mu$ and $\mu_k = k \mu$.
Pick charts $\varphi_k: \B \rightarrow X$ with $\varphi_k (0) = p_k$.
Then the curvature of the rescaled submanifold $RY_k \subset \C^n$ isn't negative with respect to the renormalized metric $(R\varphi_k) ^* \mu_k = k  (R\varphi_k) ^* \mu$.
In order to get limit submanifolds, WLOG, assume that the charts $(\varphi_k)_k$
are normal and non-degenerate (cf. Definition \ref{renorm sequence}).

Up to subsequences, the rescaled submanifolds $(RY_k)_k$ converge to some limit submanifold $Y_{\infty}$. 
For every $z\in B(\sqrt{k})$, the tangent map $d(R\varphi_k)(z)$ satisfies:
$$
\sqrt{k}\, d(R\varphi_k)(z) = d \varphi_k \left( \frac{z}{\sqrt{k}} \right).
$$
If the charts $(\varphi_k)_k$ converge to $\varphi_{\infty}$ then the linear maps $\left( \sqrt{k}\, d(R\varphi_k)(z) \right)_k$
converge to $d\varphi_{\infty} (0)$ which doesn't depend on $z$.
Hence, up to subsequences,
the renormalized Hermitian metrics converge to some constant K\"ahler metric on $\C^n$.

Since $(RY_k)$ converge to $Y_{\infty}$, the curvature of $Y_{\infty}$ at $0\in \C^n$ 
(for this constant K\"ahler metric)
is non-negative. 
This is a contradiction.

Hence, in order to prove Theorem \ref{negatively curved}, one only needs to produce a suitable subset
$A \subset \Jet^2_{d,n}$ for each type of curvature.

\subsection{Curves with negative curvature}

Here, we prove Theorem \ref{negatively curved}\ref{curve}.

\begin{prop}\label{inflection}
A complex curve $Y\subset \C^n$ has non-positive Gauss curvature
for the metric induced on $Y$ by any {\it constant} K\"ahler metric on $\C^n$.
Moreover the points where the curvature vanishes are the inflection points.
\end{prop}

\proof
Consider the Gauss map $G$ from $Y$ to $\C P^{n - 1}$. For every $y\in Y$, the value $G(y)$ is just the
tangent line $T_y Y$.
Let ${\cal O} ( - 1)$ be the usual tautological line bundle over $\C P^{n - 1}$, endowed
with the Hermitian metric induced by the constant K\"ahler metric over $\C^n$.
The tangent bundle is the pull-back $G^{*} ({\cal O}( - 1))$
of this line bundle by the Gauss map. Moreover, the metric of $TY$ equals the pull-back metric.

It is well known that the curvature of ${\cal O} ( - 1)$ is negative. Hence, since the map $G$ is holomorphic, the pull-back  $G^{*} ({\cal O}( - 1))$
has non-positive curvature. Moreover, the curvature of $G^{*} ({\cal O}( - 1))$ vanishes at critical points of $G$ ({\it i.e.} inflection points).

\qed
\bigskip

In the set $\Jet^2_{1,n}$ of $2-$jets of curves, let $A$ be the subset of $2-$jets 
of curves at inflection points.   
Notice that the affine transformations preserve $A$.
The subset
$A$ is the zero section of the vector bundle $\Jet^2_{1,n} \rightarrow \Jet^1_{1,n}$. 
Hence: $\codim (A) = n-1$. Here $d=1$ and $n\geq 3$. 
Therefore, condition $\codim (A) > d$ is satisfied and the
proof of Theorem \ref{negatively curved}\ref{curve} is complete.

\subsection{Ricci curvature}\label{proof Ricci}

Here, we prove Theorem \ref{negatively curved}\ref{item Ricci}.

\medskip

Let $Y\subset \C^n$ be a complex submanifold and denote
by $\II$  the second fundamental form of $Y$,
which depends only on the $2-$jet of $Y$.
By Proposition 9.5 in Chapter IX of \cite{KoNo69}, 
the Ricci curvature of $Y$ at a point $p$ satisfies:
\begin{eqnarray*}
\Ric (v,v) &=& -2 \sum_{i=1}^d \left\| \II(e_i, v ) \right\| ^2
\end{eqnarray*}
for every $v \in T_p Y$,
where $(e_i)_{1\leq i \leq d}$ is a unitary basis in $T_p Y$.
Hence,
if there exists $v'\in T_p Y$ such that $\II(v',v) \neq 0$,
then $\Ric (v,v)$
is negative.
This leads us to consider the set $A \subset \Jet^2_{d,n}$ of $2-$jets $a = (p_a, T_a, \II _a)$ 
such that:
$$
\exists u\neq 0 \in T_a,
\;\;\; 
\forall u'\in T_a 
\;\;\;\;\;\;\;\; 
\II _a (u,u') = 0
$$
with the notations of Section \ref{jets}.
Obviously, $A$ is invariant under the affine transformations.
We just need to show that $A$ is a closed algebraic subset of codimension $ >d$.

\medskip

For this purpose, we regard the set $A$ as the image of some subset $A_1$ of $\Jet^2_{d,n} \times \C P^{n - 1}$
which we will define below.
A point in $\Jet^2_{d,n} \times \C P^{n - 1}$ 
is a pair $(a,b)$
where $a$ is the $2-$jet of a submanifold $Y$
and $b$ is a line in $\C^n$.

The condition $b \subset T_a$ defines a closed algebraic subset $B \subset  \Jet^2_{d,n} \times \C P^{n - 1}$ 
of codimension $n - d$.
Define another closed algebraic subset $A_1 \subset B$ by using the following condition:
$$
\forall u \in b, \; \forall u'\in T_a 
\;\;\;\;\;\;\;\; 
\II _a (u,u') = 0. 
$$
In short, our condition writes: $\II _a (u,T_a) = 0$ where $u$ is a generator of the line $b$.

The map $\II _a (u,.)$ is any linear map from $T_a$ to $\C^n / T_a$.
Of course, the space of linear maps from $T_a$ to $\C^n / T_a$ has dimension $d(n - d)$.
Hence, since the subset $A_1$ is defined by $\II _a (u,.) = 0$, the codimension of $A_1$ in $B$ is $d(n - d)$
and the codimension of $A_1$ in $\Jet^2_{d,n} \times \C P^{n - 1}$  equals:
\begin{eqnarray*}
\codim (A_1) &=& d(n - d) + \codim (B) 
\\ & = & (d + 1) (n - d).
\end{eqnarray*}
Apply Theorem \ref{proj codim}. Since $A$ is the image of $A_1$ via the projection $\Jet^2_{d,n} \times \C P^{n - 1} \rightarrow \Jet^2_{d,n}$,
the subset $A$ is closed and:
\begin{eqnarray*}
\codim (A) & \geq & \codim (A_1) - \dim ( \C P^{n - 1} )
\\ & = & (d + 1) (n - d) - (n - 1)
\\ & = & d (n - d - 1) + 1.
\end{eqnarray*}
Under the assumption $d\leq n - 2$, one obtains: $\codim (A) > d$.

\subsection{Scalar curvature}

Here, we prove Theorem \ref{negatively curved}\ref{item Scal}.

\medskip

Consider the scalar curvature.
$$
\Scal = 2 \sum_{i=1}^d \Ric (e_i, e_i).
$$
Hence, the complex submanifold $Y \subset \C^n$ satisfies:
$$
\Scal = - 4 \sum_{i,j=1}^d \left\| \II (e_i, e_j) \right\| ^2.
$$
Therefore, the scalar curvature vanishes if and only if $\II$ does.
The condition $\II _a =0$ defines a codimension $\frac{d(d + 1) (n - d)}{2}$ closed algebraic subset  $A \subset \Jet^2_{d,n}$.
The assumptions $n\geq 3$ and $d\leq n - 1$ ensure: $\codim(A) > d$.

\subsection{Holomorphic sectional curvature}

Here, we prove Theorem \ref{negatively curved}\ref{item holsec}.

\medskip

The holomorphic sectional curvature HolSec$(v)$ equals $R(v,Jv,v,Jv)$ where $R$ is the usual Riemann tensor and $J$ is the multiplication by $i$ on the tangent space.
If $Y\subset \C^n$ is a complex submanifold then,
by Proposition 9.2 in Chapter IX of \cite{KoNo69}:
$$
\mbox{HolSec} (v) = -2 \| \II (v,v) \|^2.
$$
The strategy is the same as in Section \ref{proof Ricci}.
We use the same subset $B\subset  \Jet^2_{d,n} \times \C P^{n - 1}$
and some slightly different subset $A_1 \subset B$.

As before, we use the notation $(a,b)$ for a point in $\Jet^2_{d,n} \times \C P^{n - 1}$.
The subset $A_1$ is defined by the condition $\II _a(b,b) = 0$ so the codimension of $A_1$ in $B$ equals $n - d$.
Hence, by the same calculation
as in Section \ref{proof Ricci}, we get $\codim (A) \geq n - 2d + 1$ 
where $A$ is the image of $A_1$ in $ \Jet^2_{d,n}$. 
We obtain $\codim(A) > d$ because $n\geq 3d$. 

\subsection{Holomorphic bisectional curvature}\label{proof bisec}

Here, we prove Theorem \ref{negatively curved}\ref{item bisec}.

\medskip

The holomorphic bisectional curvature is defined in \cite{GoKo67}
by the formula:
$$\mbox{HolBisec}(v,v') = R(v,Jv,v',Jv').$$ 
For a submanifold $Y\subset \C^n$, formula (9) in \cite{GoKo67} writes:
\begin{eqnarray*}
\mbox{HolBisec} (v,v')  &=& - \| \II (v,v') \|^2 - \| \II (v,Jv') \|^2 
\\ &=& -2 \| \II (v,v') \|^2.
\end{eqnarray*}
We define a subset $B\subset  \Jet^2_{d,n} \times \C P^{n - 1} \times \C P^{n - 1}$ by the following characterization.
The point $(a,b,b')$ belongs to $B$ if and only if $b$, $b' \subset T_a$. 
$$
\codim (B) = 2 (n - d)
$$
Then the subset $A_1 \subset B$ is defined by the condition $\II _a(b,b') = 0$, so the codimension of  $A_1$ in $B$ is $n - d$
and the codimension of  $A_1$ in $\Jet^2_{d,n} \times \C P^{n - 1} \times \C P^{n - 1}$ equals $3(n - d)$.
As usual, $A$ is the image of $A_1$ in $ \Jet^2_{d,n}$. 
We obtain:
\begin{eqnarray*}
\codim (A) & \geq & 3(n - d) - \dim ( \C P^{n - 1} \times \C P^{n - 1} )
\\ & = & n - 3d + 2
\end{eqnarray*}
Under the assumption $n\geq 4d-1$, we obtain: $\codim (A) > d$. This completes the proof of Theorem \ref{negatively curved}.

\section{Proof of Theorem \ref{exterior}}\label{proof exterior}

\subsection{Subbundles and stationary points}

Let $X$ be a complex manifold and let $E \to X$ be a Hermitian vector bundle.
Let $\nabla^{E}$ be the Chern connection. The connection $\nabla^{E}$ is compatible with the metric and with the holomorphic structure of $E$.

Let $E_1 \subset E$ be a (smooth) complex subbundle and let $E_1^{\perp}$ be the orthogonal complement.
We write the connection $\nabla^{E}$ in the decomposition $E = E_1 \oplus E_1^{\perp}$.
$$
\nabla^{E} = \left(
\begin{array}{cc}
{\nabla}^{E_1} &  -\alpha^{*} \\ \alpha & \nabla^{E_1^{\perp}}
\end{array}
\right)
$$ 

Here, $\nabla^{E_1}$ and $\nabla^{E_1^{\perp}}$ are connections compatible with the metrics of $E_1$ and $E_1^{\perp}$.

If $E_1$ is a holomorphic subbundle, then the connection $\nabla^{E_1}$ is the Chern connection of $E_1$ 
(every holomorphic section of $E_1$ is a holomorphic section of $E$)
and, under the identification $E_1^{\perp} = E / E_1$,
the connection $\nabla^{E_1^{\perp}}$ is the Chern connection of $E / E_1$
(locally, every holomorphic section of $E / E_1$ is the projection of a holomorphic section of $E$).

Let $s$ be a section of $E_1$. A priori, the covariant derivative $\nabla^{E} s$ takes its values in $E$.
The subbundle $E_1 \subset E$ is {\it stationary} at a point $p\in X$
if, for every smooth section $s$ of $E_1$, the covariant derivative $(\nabla^{E} s)_p$ takes its values in $E_1$.
In other words, $\alpha=0$ at $p$.

Remark.
Since the connection $\nabla^{E}$ is compatible with the metric, $E_1$ is stationary at $p$
if and only if $E_1^{\perp}$ is stationary at $p$.

\medskip

At a stationary point, the curvature $F^{E_1}$  of $\nabla^{E_1}$ 
and the curvature $F^{E}$ of $\nabla^{E}$
are closely related.

\begin{lem}\label{stationary}
Let $E_1$ be a subbundle of a Hermitian bundle $E \to X$.
At a stationary point of $E_1$, every smooth sections $s_1$ and $s_2$ of $E_1$
and every vector fields $v$ and $w$ satisfy the following identity:
$$
\langle F^{E} (v,w) s_1, s_2 \rangle 
=
\langle F^{E_1} (v,w) s_1, s_2 \rangle .
$$

\end{lem}

\proof
We write the definition of the curvature.
\begin{eqnarray*}
F^{E} (v,w) s_1 &=& \nabla_v^{E}  \nabla_w^{E} s_1 - \nabla_w^{E}  \nabla_v^{E} s_1 - \nabla_{[v,w]}^{E} s_1
\\
F^{E_1} (v,w) s_1 &=& \nabla_v^{E_1}  \nabla_w^{E_1} s_1 - \nabla_w^{E_1}  \nabla_v^{E_1} s_1 - \nabla_{[v,w]}^{E_1} s_1
\end{eqnarray*}
We will compute the first summand. The section $\nabla_w^{E} s_1 - \nabla_w^{E_1} s_1$ of $E_1^{\perp}$ is orthogonal to the section $s_2$ of $E_1$:
\begin{eqnarray*}
\langle \nabla_w^{E} s_1 , s_2 \rangle & = & \langle \nabla_w^{E_1} s_1 , s_2 \rangle .
\end{eqnarray*}
The following identity holds at a stationary point:
\begin{eqnarray*}
\langle \nabla_w^{E} s_1 , \nabla_v^{E} s_2 \rangle
&=& \langle  \nabla_w^{E_1} s_1 ,  \nabla_v^{E_1} s_2 \rangle.
\end{eqnarray*}
Therefore:
\begin{eqnarray*}
\langle \nabla_v^{E}  \nabla_w^{E} s_1 , s_2 \rangle & = & d_v \langle \nabla_w^{E} s_1 , s_2 \rangle 
- \langle \nabla_w^{E} s_1 , \nabla_v^{E} s_2 \rangle
\\
& = & d_v \langle \nabla_w^{E_1} s_1 , s_2 \rangle
- \langle \nabla_w^{E_1} s_1 , \nabla_v^{E_1} s_2 \rangle
\\
& = & \langle \nabla_v^{E_1}  \nabla_w^{E_1} s_1 , s_2 \rangle.
\end{eqnarray*}
The second term is similar:
$$\langle \nabla_w^{E}  \nabla_v^{E} s_1 , s_2 \rangle
=\langle \nabla_w^{E_1}  \nabla_v^{E_1} s_1 , s_2 \rangle.$$
Since the section $\nabla_{[v,w]}^{E} s_1 - \nabla_{[v,w]}^{E_1} s_1$ of $E_1^{\perp}$ is orthogonal to the section $s_2$ of $E_1$, the third term satisfies:
$$\langle \nabla_{[v,w]}^{E} s_1 , s_2 \rangle
= \langle \nabla_{[v,w]}^{E_1} s_1 , s_2 \rangle.$$
\qed

\subsection{Griffiths-positivity and the Serre line bundle}\label{def Griffiths}

\begin{defn}
Let $X$ be a complex manifold equipped with a Hermitian vector bundle $E \rightarrow X$. 
The Hermitian bundle
 $E$ is {\it Griffiths-positive} if the curvature $F^{E}$ of $E$ satisfies
$$
 \langle F^{E} ( J v , v ) v' , J' v' \rangle > 0
$$
for all non-zero vectors $v\in T_p X$ and $v' \in E_p$, where $p$ is a point in $X$, 
where $J$ is the multiplication by $i$ on the tangent space $T_p X$
and $J'$ is the multiplication by $i$ on the fiber $E_p$ over $p$.
\end{defn}

Let $E$ be a complex vector bundle and let $P( E^{*})$ be the projectivized dual bundle.
The Serre line bundle is the line bundle ${\cal O}_{P(E^{*})} (1)$ 
on the total space of $P(E^{*})$.
A Hermitian metric on $E$ induces a Hermitian metric on ${\cal O}_{P(E^{*})} (1)$
(because ${\cal O}_{P(E^{*})} (1)$ is a quotient bundle of the pull-back of $E$ on $P(E^{*})$).
The following lemma is a reformulation of the definition of the Griffiths-positivity.
This result is known to experts.

\begin{lem}\label{SerreGriffiths}
Let $E$ be a Hermitian vector bundle and let ${\cal O}_{P(E^{*})} (1)$ be the Serre line bundle of $E$,
endowed with the corresponding metric. The Hermitian bundle $E$ is Griffiths-positive
if and only if the curvature of ${\cal O}_{P(E^{*})} (1)$ is positive.
\end{lem}

\proof 
If $E$ is Griffiths-positive, then the curvature of the pull-back of $E$ on $P(E^{*})$ is semi-positive and it is positive
in the complex directions transverse to the fibers of the fibration $P(E^{*}) \to X$.
Hence the curvature of the quotient bundle ${\cal O}_{P(E^{*})} (1)$ is positive in these directions.
Moreover, along the fibers, the curvature of ${\cal O}_{P(E^{*})} (1)$ is the usual Fubini-Study form.
In conclusion, the curvature of ${\cal O}_{P(E^{*})} (1)$ is positive in all complex directions.

Conversely, if the curvature of ${\cal O}_{P(E^{*})} (1)$ is positive, we will prove that 
$\langle F^{E} (Jv,v) s , J's  \rangle$
is positive for every non-zero vector field $v$ and every non-zero section $s$ of $E$.
In order to prove this inequality at a point $p\in X$, we will use a local holomorphic section $\varphi$ of $E^{*}$
satisfying certain conditions. 

We will choose the local section $\varphi$ to be never-zero, which can always be achieved.
The kernel of $\varphi$ is a subbundle $E_1$ of $E$.
The orthogonal complement $E_1^{\perp}$ is the pull-back $\varphi^{*} {\cal O}_{P(E^{*})} (1)$ of the Serre line bundle.
Note that the map $\varphi$ is an immersion and the curvature of ${\cal O}_{P(E^{*})} (1)$ is positive.
Hence, the curvature of the pull-back $E_1^{\perp}$ is positive.

We furthermore choose $\varphi$ such that the covariant derivative
$\nabla^{E^*} \varphi$ vanishes at $p$, which again can always be achieved.
The subbundles $E_1$ and $E_1^{\perp}$ are stationary at $p$.
If the value $s(p)$ of the section $s$ lies in fiber $(E_1^{\perp})_p$ of
the subbundle $E_1^{\perp}$ 
(for a suitable $\varphi$) then
apply Lemma \ref{stationary} to the stationary subbundle $E_1^{\perp}$. The real number
$\langle F^{E} (Jv,v) s , J's  \rangle$ is equal to the positive number $\langle F^{E_1^{\perp}} (Jv,v) s , J's  \rangle$. Therefore
the bundle $E$ is Griffiths-positive.
\qed

\medskip

Remark.
By definition, a vector bundle over a compact manifold is {\it ample} if the corresponding Serre line bundle is ample in the usual sense.
According to Kodaira's embedding theorem and Lemma \ref{SerreGriffiths}, if a Hermitian bundle $E$ is Griffiths-positive then $E$ is ample.
This result is well known and the converse is a famous conjecture.

\medskip

{\bf Conjecture} (Griffiths){\bf .}
If $E$ is an ample vector bundle then there exists a Hermitan metric $h$ on $E$
such that $(E,h)$ is Griffiths-positive.

\subsection{Griffiths-positivity of $\bigwedge ^l (TY)$}

The proof of Theorem \ref{exterior} proceeds in two steps:
(1) We study the corresponding problem in $\C ^n$.
(2) We apply the Donaldson-Auroux machinery.

\bigskip

Step 1.
Let $Y \subset \C^n$ be a complex submanifold of dimension $d$.
Consider the Gauss map $G$ defined by $G(p, [v]) = [v]$ from the total space of the projectivized bundle $P (\bigwedge ^l (TY))$
to the projective space $P(\bigwedge ^l (\C^n) )$.
Similarly, define $G'$ from the total space of $P (\bigwedge ^l (NY)^*)$
to $P(\bigwedge ^l (\C^n)^* )$.

Endow $\C^n$ with a constant K\"ahler metric. This induces Hermitian metrics on both bundles $\bigwedge ^l (T^{*} Y)$
and $\bigwedge ^l (NY)$.

\begin{lem}\label{Gauss immersive}
The Hermitian bundle $\bigwedge ^l (T^{*} Y)$ is Griffiths-positive if and only if the Gauss map $G$ is an immersion.
\end{lem}
Of course, the same holds for $\bigwedge ^l (N Y)$ and $G'$.

\proof
Denote by ${\cal O} (1)$ the usual hyperplane line bundle over $P(\bigwedge ^l (\C^n) )$. 
Endow ${\cal O} (1)$ with the Hermitian metric induced by some constant K\"ahler metric over $\C^n$.
Denote by $G^{*} ({\cal O}(1))$ the pull-back Hermitian line bundle.
Note that the Hermitian line bundle 
$G^{*} ({\cal O}(1))$ is the Serre line bundle of the Hermitian bundle $\bigwedge ^l (T^{*} Y)$.
According to Lemma \ref{SerreGriffiths},
the Hermitian bundle $\bigwedge ^l (T^{*} Y)$ is Griffiths-positive if and only if the Hermitian line bundle
$G^{*} ({\cal O}(1))$ has positive curvature.

Now, we use the fact that the curvature of ${\cal O} (1)$ is positive, as in the proof of Proposition \ref{inflection}: 
since the map $G$ is holomorphic, the pull-back  $G^{*} ({\cal O}(1))$
has semi-positive curvature and furthermore, the curvature of $G^{*} ({\cal O}(1))$ is positive if and only if $G$ is an immersion.

\qed

\bigskip

Step 2. In order to apply Theorem \ref{avoidance}, we consider some appropriate subset $A \subset \Jet^2_{d,n}$.
Note that the $1-$jet of the Gauss map of a submanifold at a point depends only on the $2-$jet of the submanifold at that point.
Let $A$ be the set of $2-$jets of submanifolds $Y\subset \C^n$ whose Gauss map $G$ is not an immersion (over some neighborhood of the corresponding point $p\in \C^n$).

\medskip

\begin{lem}\label{Gauss codim}
The set $A$ is a closed algebraic subset. Moreover, if $2d(1+l ) \leq l (n + l)$ then: $\codim (A) > d$.
\end{lem}
Similarly, if $l( n - l) \leq 2d (1 - l)$ then the codimension of the set of $2-$jets of submanifolds $Y$ with non-immersive $G'$ 
is $>d$.

\medskip

\proof
First, pick a complex submanifold $Y \subset \C^n$ of dimension $d$ and some point $p \in Y$.
The domain of $G$ is the set $P (\bigwedge ^l (TY))$ of classes of $l-$vectors in $T Y$. Denote by $\Gdec$ the restriction of $G$ to the set 
$\Pdec (\bigwedge ^l (TY))$
of classes of {\it  decomposable} $l-$vectors. The following four conditions are equivalent:

\medskip

(i) The map $G$ is an immersion in some neighborhood of the fiber $P (\bigwedge ^l (T_pY))$ in the total space $P (\bigwedge ^l (TY))$.

(ii) The restriction map $\Gdec$ is an immersion in a neighborhood of the fiber $\Pdec (\bigwedge ^l (T_pY))$ in the total space $\Pdec (\bigwedge ^l (TY))$.

(iii) For every non-zero vector $u \in T_p Y$, the subspace:
$$
K_u = \{
u' \in T_p Y,\; \II (u,u')=0
\}
$$
(where $\II$ is the second fundamental form of the submanifold $Y$ at $p$),
satisfies:  $\dim (K_u) < l$.

(iv) For every non-zero vector $u \in T_p Y$ and every linear subspace $V \subset T_p Y$ of dimension $l$, 
there exists $u'\in V$ such that $\II (u,u') \neq 0$.

\bigskip

(i) $\Rightarrow$ (ii) and equivalence (iii)$\Leftrightarrow$(iv) are clear. 

\medskip

Let's prove the contrapositive of (ii)$\Rightarrow$(iv).
Assume there exists a non-zero $u\in T_p Y$ and a linear subspace $V$ of dimension $l$, 
such that $\II (u,V) = 0$. Let $u_1$, ..., $u_l$ be a basis in $V$.
For $1\leq \alpha \leq l$,
extend $u_\alpha$ to a local section of $TY$.
The exterior product $u_1 \wedge \dots \wedge u_l$ is a non-zero decomposable
local section of $\bigwedge ^l (TY)$. Denote by $\lambda$ the induced section of  $\Pdec (\bigwedge ^l (TY))$.

Consider $u_{\alpha}$ as a map from a neighborhood of $p$ in $Y$ to $\C^n$.
The differential of $u_\alpha$ at point $p$ is a linear map
$(du_{\alpha})_p$ from $T_p Y$ to $\C^n$.

Since $\II (u,u_{\alpha})=0$, (WLOG) we can assume that we have constructed each extention $u_{\alpha}$ in such a way that
$(du_{\alpha})_p (u) = 0$. Then $d(u_1 \wedge \dots \wedge u_l)_p(u)$ is zero and $(d\lambda)_p (u)$
belongs to the kernel of the differential $(d\Gdec )_{\lambda(p)}$. 
Therefore the map $\Gdec$ isn't an immersion at point $\lambda(p)$.
This completes the proof of (the contrapositive of) (ii)$\Rightarrow$(iv).

\bigskip

In order to prove (iii)$\Rightarrow$(i), we need to look at 
the relations between the second fundamental form $\II$ and the sections of $\bigwedge ^l (TY)$.

Since $Y$ is a submanifold in $\C ^n$, a local section of $TY$ provides a map $f$ from a neighborhood of $p$ in $Y$ to $\C ^n$.
Its differential $(df)_p$ is a linear map from $T_pY$ to $\C ^n$.
If $u \in T_p Y$ is a vector then $(df)_p (u) \in \C ^n$ and, by definition of the second fundamental form,
the class of $(df)_p (u)$ mod $T_p Y$ equals $\II (u,f(p))$.

More generally, a local section of $\bigwedge ^l (TY)$ provides a local map $f$ from $Y$ to $\bigwedge ^l (\C ^n)$.
A priori, if $u \in T_p Y$ then $(df)_p (u) \in \bigwedge ^l (\C ^n)$.
Actually, $(df)_p (u)$ belongs to some linear subspace $\Gamma \subset   \bigwedge ^l (\C ^n)$.
Here, $\Gamma$ is the image of $\C ^n \otimes (T_p Y)^{\otimes l-1} \subset (\C ^n)^{\otimes l}$ 
under the usual projection $(\C ^n)^{\otimes l}  \rightarrow  \bigwedge ^l (\C ^n)$.

Let $\II \Lambda _u^l$ be the unique linear map 
from $\bigwedge^l (T_p Y)$ to $(\C^n / T_p Y) \otimes \bigwedge^{l - 1} (T_pY)$
such that:
$$
\II \Lambda _u^l (u_1 \wedge \dots \wedge u_l) = \sum_{\alpha =1}^l
(-1)^{\alpha +1} \II (u, u_{\alpha}) \otimes (u_1 \wedge \dots \wedge \widehat{u}_{\alpha}  \wedge \dots \wedge u_l)
$$
for all $u_1$, ..., $u_l\in T_pY$.
By the Leibniz rule, the class of $(df)_p (u)$ mod $\bigwedge ^l(T_p Y)$ equals $\II \Lambda_u^l(f(p))$,
where we use the following identification: $$\Gamma /  \bigwedge \hspace{-0.1cm}^l (T_p Y) =   (\C^n / T_p Y) \otimes \bigwedge \hspace{-0.1cm}^{l - 1} (T_pY).$$
The kernel of the map $\II \Lambda _u^l$ equals $\bigwedge ^l (K_u)$.

\medskip

Let's prove (iii)$\Rightarrow$(i).
Assume $(dG)_q (v) = 0$
where $q$ is a point in $P( \bigwedge ^l (T_p Y))$ and the vector $v$ is tangent to the total space $P( \bigwedge ^l (T Y))$ at $q$.
In order to show that $G$ is an immersion, we shall prove $v=0$.

Denote by $u\in T_p Y$ the projection of $v$.
The total space $P( \bigwedge ^l (T Y))$ is a subset of the product $Y \times \bigwedge ^l (\C^n)$.
The components of $v$ are $u$ and $(dG)_q (v)$. We know the second one is zero. Hence, we just need to prove $u=0$.

\medskip

Since $(dG)_q (v)$ is zero,
there exists a map $f$ from a neighborhood of $p$ to $\bigwedge ^l (\C^n)$ such that:

\medskip

\indent
- $f$ corresponds to a local section of $\bigwedge ^l (TY)$

- the vector $f(p) \in \bigwedge ^l (T_p Y)$ is a generator of the line $q \in P( \bigwedge ^l (T_p Y) )$

- $(df)_p (u) = 0$ and therefore $\II \Lambda_u^l (f(p))=0$.

\medskip

Assume $u\neq 0$. The kernel of $\II \Lambda _u^l$ equals $\bigwedge^l ( K_u )$ and assumption (iii) states that $\dim( K_u ) <l$.
Now, $f(p)\in \bigwedge ^l ( K_u )$ is a non-zero $l-$vector over a vector space of dimension $<l$.
This is a contradiction. Hence, $G$ is an immersion.
This completes equivalences (i)$\Leftrightarrow$(ii)$\Leftrightarrow$(iii)$\Leftrightarrow$(iv).

\bigskip

The equivalence (i)$\Leftrightarrow$(iv) provides another definition of the set $A$:
\begin{eqnarray*}
A &=& \{
a=(p_a,T_a, \II _a) \in \Jet^2_{d,n},\;
\exists b \in \C P^{n - 1}, \exists b' \in \mbox{Grass}(l,\C^n) \\
&  & 
\mbox{ such that }
b,b' \subset T_a
\mbox{ and }\II _a (b,b') = 0
\}.
\end{eqnarray*}
In order to study $A$, first consider some simpler sets.
The incidence subset $A_1$:
 $$
 A_1 = \{ (a,b) \in  \Jet^1_{d,n} \times \C P^{n - 1},\; b\subset T_a \} 
 $$
 has codimension $n-d$ in $\Jet^1_{d,n} \times \C P^{n - 1}$.
 Similarly:
  $$
 A_2 = \{ (a,b') \in  \Jet^1_{d,n} \times \mbox{Grass}(l,\C^n),\; b'\subset T_a \} 
 $$
 has codimension $l(n-d)$ in $\Jet^1_{d,n} \times \mbox{Grass}(l,\C^n)$ and:
$$
 A_3 = \{ (a,b,b') \in  \Jet^1_{d,n} \times \C P^{n - 1} \times \mbox{Grass}(l,\C^n),\; b,b'\subset T_a \} 
 $$
 has codimension $(1+l)(n-d)$ in $\Jet^1_{d,n} \times \C P^{n - 1} \times \mbox{Grass}(l,\C^n)$.
 
 Denote by $\pi$ the natural projection: 
 $$\pi: \Jet^2_{d,n} \times \C P^{n - 1} \times \mbox{Grass}(l,\C^n) \rightarrow \Jet^1_{d,n} \times \C P^{n - 1} \times \mbox{Grass}(l,\C^n).$$
 
 Consider:
 \begin{eqnarray*}
 A_4 &=& \{ (a,b,b') \in  \Jet^2_{d,n} \times \C P^{n - 1} \times \mbox{Grass}(l,\C^n),\; 
 \\
 & &\mbox{ such that } 
 \pi (a,b,b') \in A_3
 \mbox{ and }\II _a (b,b') = 0 \}. 
\end{eqnarray*}  
The set $A_4$ is a subbundle of the bundle $\pi: \pi^{-1}(A_3) \rightarrow A_3$ with: 
$$\dim (A_4) = \dim ( \pi^{-1}(A_3) ) - l(n - d)$$
Indeed, consider a point $(a,b,b') \in A_3 \subset \Jet^1_{d,n} \times \C P^{n - 1} \times \mbox{Grass}(l,\C^n)$.
The fiber of $\pi$ over $(a,b,b')$
equals Sym$^2(T_a, \C^n / T_a)$ and the fiber of the subbundle $A_4$ is the subspace of
the symmetric maps from $T_a \times T_a$ to $\C^n / T_a$
whose restriction to the subspace $b\times b' \subset T_a \times T_a$ is $0$.
This fiber has codimension $l(n - d)$ in Sym$^2(T_a, \C^n / T_a)$.

Hence:
$$
\codim (A_4) = \codim (A_3) + l (n - d) = (1+2l) ( n - d).
$$
The set $A$ is the image of $A_4$ by the first projection $\Jet^2_{d,n} \times \C P^{n - 1} \times \mbox{Grass}(l,\C^n) \rightarrow \Jet^2_{d,n}$.
By Theorem \ref{proj codim}, the subset $A$ is closed and satisfies:
 \begin{eqnarray*}
 \codim (A) & \geq & \codim (A_4) - \dim ( \C P^{n - 1}  \times \mbox{Grass}(l,\C^n)) \\
 & = &  (1+2l) ( n - d) - n + 1 - l (n - l) \\
 & = &  1 - d + l (n+l) - 2dl.
 \end{eqnarray*}   
Under the assumption $2d(1+l) \leq l (n+l)$,
we obtain:
$\codim (A) > d$. 

\qed

\medskip

Now we are able to complete the proof of Theorem \ref{exterior}. Our argument is almost the same as in Section \ref{preuve type}. 
Apply Theorem \ref{avoidance}, Lemma \ref{Gauss immersive} and Lemma \ref{Gauss codim}:
since $\codim (A) > d$, there exists a renormalizable sequence $(Y_k)_{k \gg 1}$ such that, 
for every limit submanifold $Y_{\infty}$, the Hermitian bundle $\bigwedge ^l (T^{*} Y_{\infty})$ is Griffiths-positive.

Since the geometry of $Y_k$ tends to the geometry of $Y_{\infty}$, it seems quite likely that, for every sufficiently large $k$, the Hermitian bundle 
$\bigwedge ^l (T^{*} Y_k)$ is Griffiths-positive. Actually, a detailed proof of this fact can be achieved by the same 
argument as in Section \ref{preuve type}: proof by contradiction and extraction of subsequences.

This completes the proof of Theorem \ref{exterior}  (of course, case (b) is similar).

\section{Proof of Theorem \ref{hyperbo}}\label{proof hyperbo}

\begin{defn}
Let $(Y_k)_{k\gg 1}$ be a renormalizable sequence of hypersurfaces.
We say the hypersurfaces  $(Y_k)_{k\gg 1}$ {\it contain lines asymptotically} if there is a limit hypersurface $Y_{\infty}$
which contains some (affine, complex) line in $\C^n$.
\end{defn}

\medskip

The proof of Theorem \ref{hyperbo} proceeds in two steps. First, by using Theorem \ref{avoidance},
we will construct a renormalizable sequence of hypersurfaces $(Y_k)_{k\gg 1}$ containing no line asymptotically.
Then, we will prove that such hypersurfaces are hyperbolic and satisfy the required estimates.

\bigskip

In order to apply Theorem \ref{avoidance}, we define some relevant closed subset in $\Jet^l_{n - 1,n}$.
Consider the set $A$ of $l-$jets of the hypersurfaces in $\C^n$ which are tangent to order $l$ 
(at least) with some affine line.

\medskip

\begin{lem}
$A$ is a closed algebraic subset in $\Jet^l_{n - 1,n}$ of codimension $ \geq l+1-n$.
\end{lem}

\proof
We will express $A$ as the projection of a subset $A_1\subset  \Jet^l_{n - 1,n} \times \C P^{n - 1}$.
Pick $(a,b) \in \Jet^l_{n - 1,n} \times \C P^{n - 1}$. 
Let $Y$ be a germ of hypersurface at a point $p$ 
with $l-$jet $a$ and let ${\cL}$ be the affine line passing through $p$
with direction $b$.
The subset $A_1$ is defined in the following way:
$(a,b)$ belongs to $A_1$ if and only if the line ${\cL}$ is tangent to order $l$ (or more) with $Y$ at $p$.
Here,  $\codim (A_1) = l$.
Clearly, $A$ is the projection  of $A_1$ in $ \Jet^l_{n - 1,n}$. 
By Theorem \ref{proj codim}, $A$ is a closed algebraic subset with:
\begin{eqnarray*}
\codim (A) & \geq & l - \dim (\C P^{n - 1})
\\ & = & l +1 - n.
\end{eqnarray*}

\qed

\begin{prop}\label{no line}
Let $X$ be a projective manifold and let $L$ be an ample line bundle over $X$.

Then, for every sufficiently large $k$, there exists a hypersurface $Y_k$
of degree $k$ 
such that the hypersurfaces $(Y_k)_{k\gg 1}$ form a renormalizable sequence and contain no line asymptotically.
\end{prop}

\proof
Notice that $A$ is invariant under the affine transformations.
Pick some large $l$, so that $\codim (A) > d = n - 1$. By Theorem \ref{avoidance},
there exist renormalizable hypersurfaces $(Y_k)_{k\gg 1}$ whose limit hypersurfaces 
have no contact of order $l$ (or more) with any line.
In particular, they contain no line. 
\qed  

\bigskip

Brody has proved a clever reparametrization lemma:

\begin{lem}\label{LemBro}
For every holomorphic map $f: \D \rightarrow X$, there exists a holomorphic map $g: \D \rightarrow X$ such that:
\begin{eqnarray*}
g(\D) & \subset & f(\D)
\\
\sup_{\D} \| g' \| & \leq & C \| g'(0) \|
\\
\| f'(0) \| & \leq & C \| g'(0) \|
\end{eqnarray*}
where $C$ is a universal bound (in particular, independent of $f$).
\end{lem}  
\proof
The idea of the proof is to go back and forth between the Euclidean metric and the Poincar\'e metric.
First, remark that the Jacobian of the map $p\mapsto \frac{p}{2}$ from $\D$ endowed with the Poincar\'e metric
to $\D$ endowed with the Euclidean metric is a positive function on $\D$ which tends to zero
on the boundary $\partial \D$.   
Therefore the jacobian $j_1$ of the map $f_1(p) = f\left( \frac{p}{2} \right)$ from $\D$ endowed with the Poincar\'e metric
to $X$ tends to zero and, hence, the supremum of $j_1$ is achieved at some point $p_0$ in $\D$.
Moreover, $\| f'(0) \| = C_1 j_1 (0)$ where $C_1$ is some universal constant.

Pick a conformal transformation $h: \D \rightarrow \D$ with $h(0)=p_0$. We set $f_2 = f_1 \circ h$. 
Since conformal mappings preserve 
the Poincar\'e metric, 
the corresponding Jacobian $j_2$ satifies $j_2 = j_1 \circ h$, so that:
$$
\sup_{\D} j_2 = \sup_{\D} j_1 = j_1(p_0) = j_2 (0). 
$$

To complete the proof, we go back to the Euclidean metric.
The Jacobian of the map $p\mapsto \frac{p}{2}$ from $\D$ endowed with the Euclidean metric
to $\D$ endowed with the Poincar\'e metric is bounded by a constant $C_2$.
Hence the map $f_3(p) =  f_2\left( \frac{p}{2} \right)$ satisfies:
$$
\sup_{\D} \| f_3' \|  \leq  C_2 \, \sup_{\D} j_2 = C_2 j_2 (0) = C_2 C_3 \| f_3' (0) \|
$$
and:
$$
\| f'(0) \| = C_1 j_1(0) \leq C_1 j_1 (p_0) = C_1 j_2(0) = C_1 C_3  \| f_3' (0) \|
$$
where $C_3$ is another constant.  
The proof is completed (set $g=f_3$). 

\qed  

\medskip

Theorem \ref{hyperbo} instantly follows from Proposition \ref{no line} and from the following renormalized version of Brody's techniques:
\begin{prop}\label{Brody}
Let $(Y_k)_{k\gg 1}$ be renormalizable hypersurfaces. 
Assume that for infinitely many $k$, 
there exists a holomorphic map $f_k: \D \rightarrow Y_k$ such
that the numbers $\sqrt{k} \| f_k'(0) \|$ tend to $\infty$.

Then the hypersurfaces $(Y_k)_{k\gg 1}$ contain lines asymptotically.
\end{prop}  
\proof
Using Lemma \ref{LemBro}, we can assume (WLOG) that:
$$
\sup_{\D} \| f_k' \| \leq C \| f_k'(0) \|
$$
for a constant $C$ which doesn't depend on $k$. 
By the assumption on $\sqrt{k} \| f_k'(0) \|$, the vector $f_k' (0)$ is large with respect to the rescaled metric $\mu_k = k \mu$.
We can also assume that $f_k' (0)$ is small with respect to the metric $\mu$.
Indeed, if $\| f_k'(0) \|$ doesn't tend to $0$, and hence
$ \| f_k'(0) \| \geq k^{-\frac{1}{4}}$ for infinitely many $k$, 
then replace $f_k$ with the rescaled map 
 $z \mapsto f_k\left( k^{-\frac{1}{4}} \| f_k'(0) \| ^{-1} z \right)$ .
The diameter of $f_k (\D)$ 
(with respect to $\mu$)
tends to $0$.

Pick some holomorphic chart $\varphi_k: \B \rightarrow X$ centered at $f_k(0)$.
We can assume that the family $(\varphi_k)$ is normal and non-degenerate. Moreover, we can assume (WLOG) that every limit map $\varphi_{\infty} : \B \rightarrow X$ 
of the normal family is a chart.
Then, the restrictions $\varphi_k |_{ \B \left( \frac{1}{2} \right) }$ 
from the ball $\B \left( \frac{1}{2} \right)$
with radius $\frac{1}{2}$ to $\varphi_k \left( \B \left( \frac{1}{2} \right) \right)$
have bounded distortion. 
The non-degeneracy of the family $(\varphi_k )$
implies that $\varphi_k \left( \B \left( \frac{1}{2} \right) \right)$ contains a ball 
of radius independent of $k$ centered at $f_k (0)$ in $X$.
Hence, for $k$ sufficiently large:
$$
f_k (\D) \subset  \varphi_k \left( \B \left( \frac{1}{2} \right) \right)
$$
Define $g_k: \D \rightarrow \varphi_k^{-1} (Y_k) \subset \C^n$ by the formula:
$$
g_k = \varphi_k^{-1} \circ f_k
$$

Since $\varphi_k$ has bounded distortion
on $\B \left( \frac{1}{2} \right)$, 
the maps $f_k$ and $g_k$ behave similarly, in particular:
\begin{enumerate}
\item $\sqrt{k} \| g_k'(0) \|$ tends to $\infty$
\item $\sup_{\D} \| g_k' \| \leq C \| g_k'(0) \|$ for some $C$.
\end{enumerate}

We now define a rescaled version $h_k$ of $g_k$. We set $r_k = \sqrt{k} \| g_k'(0) \|$
and we define $h_k$ on a large disk $\D (r_k)$ by the following formula:
$$
h_k(z) = \sqrt{k} \, g_k\left( \frac{z}{r_k} \right) \in \C^n
$$

Now:
\begin{eqnarray*}
\| h_k'(0) \| & = & \frac{\sqrt{k}}{r_k} \| g_k' (0) \| = 1
\\
\sup_{\D (r_k)} \| h_k' \| & = & \frac{\sqrt{k}}{r_k} \sup_{\D} \| g_k' \|
\\ & \leq & C \frac{\sqrt{k}}{r_k} \| g_k' (0) \| = C
\end{eqnarray*}

We conclude that the maps $(h_k)$ form a normal family because $\| h_k' \|$
is uniformly bounded and $h_k (0)=0$. Therefore some subsequence converges to a map $h_{\infty}$
with bounded derivative and $\| h_{\infty}' (0) \| =1$. 
Since $r_k$ tends to $\infty$, 
 $h_{\infty}$ is an entire map.

\medskip

Recall a classical fact. The image of a non-constant entire map 
$h_{\infty}: \C \rightarrow \C^n$ with bounded derivative is a line 
(since the derivative is bounded, Liouville's theorem implies it is constant).

The image $g_k(\D)$ is contained in $\varphi_k^{-1} (Y_k)$,
so $h_k ( \D (r_k))$ is contained in $RY_k$. Passing to a subsequence if necessary, the rescaled 
hypersurfaces $(RY_k)$ converge to some limit hypersurface $Y_{\infty}$
and the line $h_{\infty} (\C )$ is contained in $Y_{\infty}$. 
This completes the proof of Proposition \ref{Brody} (and of Theorem \ref{hyperbo}).

\qed

\section{Proof of Theorem \ref{avoidance}}\label{la theorie}

We will give the proof of Theorem \ref{avoidance}. The steps of our proof are the same as in \cite{Do96}. 
We will also use the argument proposed in \cite{Au02} to simplify a difficult step of the proof. 

\subsection{Framed charts}

Let $X$ be a projective manifold of dimension $n$ and let $L \to X$ be an ample line bundle endowed with a Hermitian metric
of curvature $-i2\pi \omega$ where $\omega$ is a K\"ahler form. We define a Riemannian metric $g$ by
the usual formula:
$$
g(v,w) = \omega (v, Jw).
$$
Hence, the tangent space $T_p X$ is a Hermitian vector space.

\medskip

Let  $\B (r) \subset \C^n$ denote the ball of radius $r$. 

\begin{defn}
A {\it framed chart} is a pair $(\varphi, \tau)$
where $\varphi: \B (r) \to X$ is a holomorphic chart
and $\tau: \C \times \B (r) \to \varphi ^{*} L$ is a holomorphic trivialization of the pull-back line bundle $\varphi ^{*} L$.
\end{defn}
We denote by $\tau^k$ and $m\tau^k$ the corresponding trivializations of the
line bundle  $\varphi ^{*} L ^{\otimes k}$ and of the vector bundle $\varphi ^{*} ( \C ^m \otimes L ^{\otimes k})$.

Let $s_k$ be a holomorphic (or smooth)
section of  $ \C ^m \otimes L ^{\otimes k}$.
Via the trivialization $m\tau^k$, the pull-back section $\varphi ^{*} s_k$ of  $\varphi ^{*} ( \C ^m \otimes L ^{\otimes k})$
defines a holomorphic (resp. smooth)
map $M(s_k, \tau): \B(r) \to \C^m$.
Then, the renormalized map $R(s_k, \tau): \B(\sqrt{k} r) \to \C^m$ 
if defined as the composite map of $M(s_k, \tau)$ with the homothety $\B(\sqrt{k} r) \to \B(r)$
of {\it ratio} $\frac{1}{\sqrt{k}}$.
We use shorter notations $Ms_k$ and $Rs_k$ instead of $M(s_k, \tau)$ and $R(s_k, \tau) $ whenever $\tau$ is clear from the context.
Hence:
$$
R s_k (v) = M s_k \left( \frac{v}{\sqrt{k}} \right)
$$
for all $v \in \B (\sqrt{k} r)$.
\medskip

Recall that $L$ is endowed with a Hermitian metric. Denote by $e^{-h}$ the metric on the trivial bundle $\C \times \B (r) \to \B (r)$
pull-back of the metric of $L$ via the trivialization $\tau$.
This defines a smooth map $h: \B (r) \to \R$.

\begin{defn}
The framed chart $(\varphi, \tau)$ is {\it standard at the origin} if the following two conditions hold:
\begin{enumerate}[label=(\alph*)]

\item
The tangent map $d \varphi (0): \C^n \to T_{\varphi (0)} X$
is an isometry of Hermitian vector spaces, where $T_{\varphi (0)} X$ is endowed with the metric $g_{\varphi (0)}$
and $\C^n$ with the standard metric.

\item
The function  $h$ has a non-degenerate maximum at $0$ with value $h (0) = 0$
and the second derivative $d^2 h(0)$ satisfies:
$d^2 h (0) (v,w) =  \pi \langle v,w \rangle $
for all vectors $v, w \in \C^n$
where $\langle.,.\rangle: \C^n \times \C^n \to \R$ denotes the usual inner product.
\end{enumerate}
\end{defn}
Notice that condition (b) implies that: $h(v) =  \frac{\pi}{2} \| v \| ^2 + o (  \| v \| ^2  )$ and hence:
$$
 \frac{\pi}{4} \| v \| ^2 \leq h(v) \leq  \frac{3\pi}{4} \| v \| ^2
$$
on a neighborhood of the origin.

\medskip

For every $p\in X$, there exists a framed chart  $(\varphi_p, \tau_p)$ centered at $p$,
standard at the origin, because the curvature of $L$ equals $-i2\pi \omega$.

\begin{defn}
For every index $i$ in a set $I$,
let $\varphi_i: \B(r) \to X$ be a holomorphic chart defined on a ball $\B(r)$ independent of $i$.
The family $(\varphi_i)_{i\in I}$ is {\it normal} if 
the charts $(\varphi_i)_{i\in I}$ form a  relatively compact subset in the set of maps $\B(r) \to X$
for the smooth compact-open topology.

Similarly, let $(\varphi_i, \tau_i)$ be framed charts defined on a ball $\B(r)$.
For every $i\in I$, the map $\tau_i: \C \times \B(r) \to \varphi_i^{*}L$
induces a map $\widetilde{\tau}_i:\C \times \B(r) \to L$
where we denote by $L$ the total space of the line bundle.
The family of framed charts $(\varphi_i, \tau_i)_{i\in I}$ is {\it normal} if the maps  $(\widetilde{\tau}_i)_{i\in I}$
form a relatively compact family in the set of maps $\C \times \B(r) \to L$
for the smooth compact-open topology.
\end{defn}

\subsection{Renormalizable sequences of sections}

In Section \ref{limit sub}, we have defined renormalizable sequences of submanifolds.
Similarly, we will define renormalizable sequences of sections.

\begin{defn}
For every sufficiently large $k$, let $s_k$ be a holomorphic section of  the vector bundle $\C ^m \otimes L ^{\otimes k}$ where the rank $m$ is independent of $k$.

The sequence $(s_k)_{k\gg 1}$ is {\it renormalizable} if for all normal sequence
of framed charts $(\varphi_k , \tau_k )_{k \in I}$ standard at the origin  (defined on a ball $\B (r)$ and 
indexed by some infinite set $I$ of sufficiently large integers), the corresponding family 
$(Rs_k: \B (\sqrt{k} r) \to \C^m)_{k \in I}$
of rescaled maps is normal.

\medskip

The limits of (subsequences of) the normal family $(Rs_k : \B (\sqrt{k} r) \to \C^m)_{k \in I}$ are holomorphic maps $s_{\infty}: \C^n \to \C^m$.
Such a $s_{\infty}$ is called a {\it limit section}.
\end{defn}

Remark. The map $s_{\infty}: \C^n \to \C^m$ is a section of the trivial bundle $\C^m \times \C^n \to \C^n$. This trivial bundle should be viewed as an "infinitesimal bundle".
Though this bundle is trivial, it is naturally endowed it with the non-trivial Hermitian metric $e^{-h_{\infty}}$ where the map $h_{\infty}: \C^n \to \R$ is defined by the formula
$h_{\infty} (v) =  \frac{\pi}{2} \| v \| ^2$. Hence, the norm of a constant section $s$ of $\C^m \times \C^n \to \C^n$
is a Gaussian function $ \| s(v) \| = \| s(0) \| \exp \left( -{ \frac{\pi}{2} \| v \| ^2} \right) $.  

\bigskip
A sequence $(s_k)_{k\gg 1}$ of sections is renormalizable if and only if, for every $l$, the corresponding sequence $(\|s_k\| _{{\cal C}^{l} })_{k\gg 1}$ of norms 
${\cal C}^{l}$ measured with the rescaled metric $g_k = k g$ is bounded. 
Of course, the definition of the norm ${\cal C}^{l}$ for sections uses connections. We use the Levi-Civita connection over $X$ 
and the Chern connection on $\C^m \otimes L^{\otimes k}$. Notice that both metrics $g$ and $g_k$ induce the same Levi-Civita connection.

\subsection{Avoidance theorem for sections}

The jet of order $l$ at the origin of a (germ of) holomorphic map $\C^n \to \C^m$
is a polynomial map  $\C^n \to \C^m$ of degree $\leq l$.
Hence, the space $\Pol _l (\C^n,\C^m)$ of polynomial maps $\C^n \to \C^m$ of degree $\leq l$
equals the space of jets at the origin.
Now, we will consider jets at any point in $\C^n$.
The space of jets of order $l$ of holomorphic maps $\C^n \to \C^m$ equals $\C^n \times \Pol _l (\C^n,\C^m)$.

\begin{thm}\label{avoidance sections}
Let $X$ be a projective manifold of dimension $n$
and let $L \to X$ be an ample line bundle endowed with a Hermitian metric
of positive curvature.
Let $A \subset  \C^n \times \Pol _l (\C^n,\C^m)$ be a closed complex algebraic subset 
of codimension $\geq n + 1$ satisfying the following invariance conditions:

\begin{enumerate}

\item
$A$ is invariant under the natural action of affine unitary transformations of $\C^n$
upon $\C^n \times \Pol _l (\C^n,\C^m)$.

\item
Let $f: (\C^n , p) \to \C^m$ be the $l-$jet of a holomorphic map and let $u: (\C^n , p) \to \C$ be the $l-$jet of a holomorphic function.
If $f$ belongs to $A$ then $uf$ belongs to $A$.

\end{enumerate}

Then for every sufficiently large $k$, there exists a holomorphic section $s_k$ of $\C^m \otimes L^{\otimes k}$ 
such that
the sections $(s_k)_{k\gg 1}$ form a renormalizable sequence 
which {\rm avoids $A$ asymptotically} in the sense that  
the jets of order $l$ of its limit sections lie in the complement $ \left( \C^n \times \Pol _l (\C^n,\C^m)  \right)  \backslash A$.

\end{thm}

This result has similarities with the main theorem in \cite{Au01}.
We will give the proof of Theorem \ref{avoidance sections}
in section \ref{do co}. Then, in section \ref{proof sub}, we will prove that this theorem implies Theorem \ref{avoidance}.

\subsection{Donaldson's construction}\label{do co}

Here, we prove Theorem \ref{avoidance sections}.

\medskip

First, we will define tools that we will use in order to prove
Theorem \ref{avoidance sections}. 
For every $p\in X$, there exist a 
framed chart $(\varphi_p, \tau_p)$ centered at $p$ and standard at the origin.
Moreover,  
if we denote by $e^{-h}$ the pull-back Hermitan metric on the trivial bundle $\C \times \B (r) \to \B (r)$,
we can assume that the map $h: \B (r) \to \R$
satisfies the condition  $\frac{\pi}{4} \| v \| ^2 \leq h(v) \leq  \frac{3\pi}{4} \| v \| ^2$
on $\B (r)$.

Since $X$ is compact, we can assume that the charts $( \varphi_p: \B (r) \to X )_{p \in X}$ are defined on a ball $\B (r) $ 
whose radius does not depend on $p$ and we can assume that the framed charts 
$(\varphi_p, \tau_p)_{p \in X}$ form a normal family.

If $s_k$ is a holomorphic section of  $\C ^m \otimes L ^{\otimes k}$, we denote by $R_p s_k$ 
the corresponding renormalized holomorphic map $\B (\sqrt{k} r) \to \C^m$ with respect to the chart $\varphi_p$ and 
the trivialization $\tau_p$.
The jet of order $l$ at the origin of $R_p s_k$ is a polynomial map $H: \C^n \to \C^m$ of degree $\leq l$.

We fix a norm $\| . \|$ on $\Pol _l (\C^n,\C^m)$ (they are all equivalent).
The Riemannian metric $g$ induces a renormalized metric $g_k = k g$ and we denote by $d$ and $d_k = \sqrt{k} d$ the corresponding distance functions. 

\begin{lem}

For every sufficiently large $k$ and for every $H \in \Pol _l (\C^n,\C^m)$, there exists a global holomorphic section $\sigma_k = \sigma_k (H,p)$
satisfying the following conditions: 

\begin{enumerate}
\item
The jet of order $l$ of $R_p \sigma_k$ at the origin equals $H$.

\item
The mass of $\sigma_k$ is concentrated around $p$, in the sense that,
for every $l'$,
the norm of the jet of order $l'$ of $\sigma_k$ at $q$ is bounded by:
$$
C \| H \|  \exp \left( - \frac{k d(p,q)^2}{C} \right) 
$$
where the constant $C$ only depends on $l$ and $l'$ (in particular, $C$ is independent of $k$ and $H$).
\end{enumerate}

\end{lem}

\proof
Let $\widetilde{\sigma}_k = \widetilde{\sigma}_k (H, p)$ be the smooth section of $\C^m \times L^{\otimes k}$ defined by the following two conditions:
\begin{enumerate}

\item
The rescaled smooth map $R_p \widetilde{\sigma}_k: \B ( \sqrt{k} r) \to \C^m$ satifies:
$$
R_p \widetilde{\sigma}_k (v) = \rho \left( \frac{v}{\sqrt{k}}  \right) H (v) 
$$
where $\rho: \B (r) \to \R_{\geq 0}$ is a cut-off function.

\item
We extend $\widetilde{\sigma}_k$ by $0$ on the complement $X \backslash \varphi (\B (r)) $.

\end{enumerate}

Recall that the pull-back metric $e^{-h}$ of the metric of $L$ satisfies $ \frac{\pi}{4} \| v \| ^2 \leq h(v) \leq  \frac{3\pi}{4} \| v \| ^2$
on  $\B(r)$.
The pull-back of the metric of  $L^{\otimes k}$ is $e^{-k h}$ and the
pull-back of the metric of $\C^m \otimes L^{\otimes k}$ equals $e^{-k h} \langle.,. \rangle$ where $ \langle.,. \rangle$ is the usual inner product of $\C^m$.
Hence, by the Leibniz rule, the norm of the jet of order $l'$ of $\widetilde{\sigma}_k$ at $q$ is bounded by $C \| H \| \exp \left( - \frac{k d(p,q)^2}{C} \right) $.
Here, we measure the norms of jets with the metric $g_k$.

The support $\supp \left( \overline{\partial}  \widetilde{\sigma}_k \right)$ is contained in the image $\varphi ( \supp ( d \rho))$.
Hence, every $q \in \supp \left(\overline{\partial}  \widetilde{\sigma}_k \right)$ satisfies $d(p,q) > \frac{1}{C}$ for a constant $C>0$ and 
the norm of the jet of order $l'$ of $\overline{\partial}  \widetilde{\sigma}_k$ at $q$ is bounded by $C \| H \| \exp \left( - \frac{k}{C} \right) $ for a (new) constant $C$.

We solve the Cauchy-Riemann equation $ \overline{\partial}  \xi_k = -  \overline{\partial}  \widetilde{\sigma}_k$. 
There exists a section $\xi_k = \xi_k (H, p) $ such that $\widetilde{\sigma}_k +  \xi_k$ is holomorphic and:
$$
\| \xi_k \|_{L^2 (X, g_k)} \leq C' \| \overline{\partial}  \widetilde{\sigma}_k  \|_{L^2 (X, g_k)} \leq C'' \| H \| \exp \left( - \frac{k}{C} \right)
$$
(estimate (35) in \cite{Do96}). A priori, the result from \cite{Do96} holds for sections of the line bundle $L^{\otimes k}$ but, clearly,  
we can extend this result to sections of the direct sum $\C^m \otimes L^{\otimes k}$. 
We can assume that $\xi_k (H, p)$ depends linearly on $H$ (choose values on a basis of $ \Pol _l (\C^n,\C^m)$ and take linear combinations).
We set $ \overline{\sigma}_k  = \widetilde{\sigma}_k +  \xi_k $. 
When $k$ tends to $\infty$, the norm $L^2$ of $\xi_k$ and the norms ${\cal C}^{l'}$ of $\overline{\partial} \xi_k$ decrease exponentially.
Recall the following elliptic estimate: 
$$
\| \xi_k \|_{ {\cal C}^{l'} (X, g_k)} \leq  C \left( \| \xi_k \|_{ L^2 (X, g_k)}  + \| \overline{\partial} \xi_k \|_{ {\cal C}^{l'} (X, g_k)}  \right) .
$$
Hence, the norm ${\cal C}^{l'}$ of $\xi_k$ satisfies a similar estimate:
$$
\| \xi_k \|_{ {\cal C}^{l'} (X, g_k)} \leq C \| H \| \exp \left( - \frac{k}{C} \right) .
$$
Collecting the estimates of $\widetilde{\sigma}_k$ and  $\xi_k$,
we obtain the following bound for the norm of the jet of $\overline{\sigma}_k  = \widetilde{\sigma}_k +  \xi_k$ at $q$:
$$
C \| H \|  \exp \left ( - \frac{k d(p,q)^2}{C} \right) 
$$
where $C$ is another constant.

Because the ${\cal C}^{l}-$norm of $\xi_k$ is bounded by $C \| H \|  \exp \left ( - \frac{k}{C} \right)$,
the jet of the renormalized holomorphic map
$R_p  \overline{\sigma}_k (H,p)$ at the origin is closed but not equal to $H$.
Hence, if we denote by $I_k(H)$ the jet of $R_p  \overline{\sigma}_k (H,p)$ at the origin,
the linear map $I_k: \Pol _l (\C^n,\C^m) \to  \Pol _l (\C^n,\C^m)$ lies
in a small neighborhood (of radius $C  \exp \left ( - \frac{k}{C} \right)$)
of the identity map if $k$ is sufficiently large. In particular $I_k$ is an isomorphism of bounded distortion and we can set: 
$$\sigma_k (H, p) =  \overline{\sigma}_k (I_k^{-1} (H), p) .$$

\qed

The mass of the section $\sigma_k (H, p)$ is concentrated around the center of the chart $\varphi_p$.
The main tools in Donaldson's construction are suitable combinations of such concentrated sections.

For every $k$, let $X_k \subset X$ be a $g_k -$discretization of $X$, 
that is, a finite subset satisfying the following two conditions:

\begin{enumerate} 

\item
For every $p \in X$, there exists $q \in X_k$ such that $d_k (p,q) < 1$

\item
Every $p$ and $q \in X_k$ satisfy  $d_k (p,q) \geq 1$

\end{enumerate}
where $d_k = \sqrt{k} d$ is the distance function of the rescaled metric $g_k$.

\begin{lem}\label{somme de pics}
For every sufficiently large $k$ and for every $p\in X_k$, let $H_{k,p} \in  \Pol _l (\C^n,\C^m)$
be a polynomial map such that $ \| H_{k,p} \| \leq 1$.
We define a section $s_k$ of $\C^m \otimes L^{\otimes k}$ by the following formula:
$$
s_k = \sum_{p \in X_k}  \sigma_k (H_{k,p}, p).
$$

Then the sections $(s_k)_{k \gg 1}$
form a renormalizable sequence.
In particular, if $K \subset \C^n$ is a compact subset
then, for every $q\in X$ and $z \in K$, the norm of the jet of order $l$  of $R_q s_k$ at $z$
is bounded by a constant independent of $k$, $q$ and $z$.
\end{lem}

\proof
The norm of the jet of order $l'$ of $\sigma_k (H_{k,p}, p) $ at $q$ is bounded by
$
C \exp \left ( - \frac{ d_k(p,q)^2}{C} \right)
$.
Hence the norm of the jet of $s_k$ is bounded by:
$$ C \sum_{p \in X_k}  \exp \left ( - \frac{d_k (p,q)^2}{C} \right).
$$

We fix $q \in X$ and we define a partition:
$$
X_k = \bigcup_{a \geq 0} X_{k,a}
$$
where, for every integer $a$,
the subset $X_{k,a}$ is the set of points $p\in X_k$ satisfying $a \leq d_k (p,q) < a+1 $.
By a volume argument, the cardinality of $X_{k,a}$ is bounded by a polynomial $P (a)$
where $P$ is independent of $k$ and $q$.
The norm of the jet of $s_k$ is bounded by:
$$
 C \sum_{a \geq 0} \; \sum_{ p \in X_{k,a}}  \exp \left ( - \frac{ d_k (p,q)^2}{C} \right)
\leq 
 C \sum_{a \geq 0}  \exp \left ( - \frac{a^2}{C} \right) P (a).
$$

This bound is finite and independent of $k$.

\qed

\medskip

Remark.
If $H_{k,p} = 0$ for every $p$ such that $d_k (p,q) \leq D$,
then the norm of the $l-$jet of $R_q s_k$ at $z$ is bounded by:
$$
 C \sum_{a \geq {D - 1}}  \exp \left ( - \frac{a^2}{C} \right) P (a)
$$
and hence by:
$$
C \exp \left( - \frac{D^2}{C} \right)
$$
where $C$ is another constant.

\bigskip

As we have seen,  the space of jets of order $l$ of holomorphic maps $\C^n \to \C^m$ equals $\C^n \times \Pol _l (\C^n,\C^m)$.
We define the (fibered) distance between two subsets $A$ and $B$ of $\C^n \times \Pol _l (\C^n,\C^m)$. 
By definition, the distance $\dist (A,B)$ is the {\it infimum}
of the distances between pairs of points $(z,H_1) \in A$ and $ (z,H_2)\in B$ having the same projection $z \in \C^n$.

\begin{lem}\label{local avoidance}
Let $A \subset  \C^n \times \Pol _l (\C^n,\C^m)$ be a closed complex algebraic subset 
of codimension $\geq n + 1$ and let $K \subset \C^n$ be a compact subset.
For every sufficiently large $k$ and for every $p\in X_k$, let $H_{k,p} \in  \Pol _l (\C^n,\C^m)$
be a polynomial map such that $ \| H_{k,p} \| \leq 1$.
We set:
$$
s_k = \sum_{p \in X_k}  \sigma_k (H_{k,p}, p).
$$

Then, 
for every $q\in X$ and $0< \varepsilon < \frac{1}{4}$, there exists a polynomial map $H \in \Pol _l (\C^n,\C^m)$
satisfying $ \| H \| \leq \varepsilon$ such that, setting
$$
\widetilde{s}_k = s_k + \sigma_k (H,q),
$$
the distance between $A$ and the set $B$ of jets of order $l$ of $R_q \widetilde{s}_k$
at points in $K$
satisfies the following estimate:
$$\dist (A,B) \geq \varepsilon ( - \log \varepsilon )^{-N}$$ 
where $N$ is independent of $k$ and $q$ and of the polynomial maps $H_{k,p}$. 
\end{lem}

\proof 
We can assume (WLOG) that $K$ is a polydisk. 
We approximate the holomorphic map $f = R_q s_k: \B (\sqrt{k} r) \to \C^m$ by its polynomial expansion $\overline{f}$.
More precisely, for every $0 < \eta < \frac{1}{4}$, there exists a polynomial map $\overline{f}: \C^n \to \C^m$
satisfying the following two conditions:

\begin{enumerate}

\item
For every $z \in K$, the Taylor polynomial $T_l (f - \overline{f}, z)$ 
of order $l$ of the map $f-\overline{f}$ at $z$ satisfies the following estimate:
$$ \left\|  T_l (f - \overline{f}, z)  \right\| \leq \eta .$$

\item
The degree of $\overline{f}$ is $\leq ( - \log \eta)^N$ where $N$ is independent of $k$, $q$ and of the maps $H_{k,p}$. 

\end{enumerate} 

For every $z \in \C^n$ and every $H \in \Pol _l (\C^n,\C^m)$, we set: $$\psi (z,H) = H - T_l (\overline{f} , z).$$
The map $\psi:   \C^n \times \Pol _l (\C^n,\C^m) \to \Pol _l (\C^n,\C^m)$ is a polynomial map of degree $\leq ( - \log \eta)^N$.

Notice that, since:
$$
\dim (A) \leq \dim ( \C^n \times \Pol _l (\C^n,\C^m) ) - (n + 1) < \dim ( \Pol _l (\C^n,\C^m) ),
$$
the restriction map $\psi_A$ isn't surjective and $\psi (A)$ is contained in a hypersurface ${\cal H} \subset \Pol _l (\C^n,\C^m)$
whose degree is a polynomial function of the degree of $\psi$, say $( - \log \eta)^N$ where $N$ is another constant.

\medskip

We denote by ${\cal B} (\varepsilon) \subset \Pol _l (\C^n,\C^m)$ the ball of radius $\varepsilon$.
The $2\eta-$neighborhood ${\cal H}_{2\eta}$ of ${\cal H}$ satisfies the following estimate:
\begin{eqnarray*}
\vol \left( {\cal H}_{2\eta} \cap {\cal B} (\varepsilon)  \right) 
& \leq & C \eta \varepsilon^{2\delta - 1} ( \deg {\cal H})^C
\\
& \leq & C \eta \varepsilon^{2\delta - 1} ( - \log \eta )^{CN}
\end{eqnarray*}
where $C$ is a constant and $\delta$ is the complex dimension of $\Pol _l (\C^n,\C^m)$.

By a simple volume argument, if we set $\eta = \varepsilon ( - \log \varepsilon )^{-N}$
where $N$ is sufficiently large then there exists $H \in \Pol _l (\C^n,\C^m)$ such that $\| H \| < \varepsilon$ and $H \notin {\cal H}_{2\eta}$. 

We set $\widetilde{s}_k = s_k + \sigma_k (H,q)$. For every $z$:
$$
T_l (R_q \widetilde{s}_k , z ) =  H + T_l (R_q s_k , z ) = H + T_l ( f , z ). 
$$

Then, for every $(z, H_1) \in A$:
\begin{eqnarray*}
\left \|  T_l (R_q \widetilde{s}_k , z )  - H_1\right\| 
& = & 
\left \|  H + T_l (f , z ) - H_1  \right\| 
\\
& \geq &
\left \|  H - \psi (z,  H_1)  \right\|  - \left \|   T_l (f - \overline{f}, z )   \right\| 
\\
& \geq &
2 \eta - \eta = \eta = \varepsilon ( - \log \varepsilon )^{-N} .
\end{eqnarray*}

\medskip

Moreover, we can assume that $N$ is independent of the polynomial maps $H_{k,p}$.
Indeed, if we consider all maps $s_k$, for all sufficiently large $k$ and for all polynomial maps $(H_{k,p})_{p\in X_k}$
such that $ \| H_{k,p} \| \leq 1$, then the corresponding maps $(R_q s_k)$ form a normal family and so we can assume that
every constant in this proof is independent of the maps $H_{k,p}$.

\qed

\begin{lem}\label{globalisation}
Let $A \subset  \C^n \times \Pol _l (\C^n,\C^m)$ be a closed complex algebraic subset 
of codimension $\geq n + 1$ and let $K \subset \C^n$ be a compact subset.

Then for every sufficiently large $k$, there exists a holomorphic section $s_k$ of $\C^m \otimes L^{\otimes k}$ 
such that
the sections $(s_k)_{k\gg k}$ form a renormalizable sequence and, 
for every $q\in X_k$, 
the distance between $A$ and the set $B$ of jets of order $l$ of $R_q s_k$
at points in $K$
satisfies $\dist (A,B) \geq \varepsilon$ where $\varepsilon >0$ is independent of $k$ and $q$. 
\end{lem}

\proof 
We will construct a section:
$$
s_k = \sum_{p\in X_k} \sigma_k (H_{k,p} , p)
$$
where the maps $H_{k,p}$ are suitable elements in $ \Pol _l (\C^n,\C^m)$. 
The construction follows the lines of \cite{Do96} and we will skip some calculations.

For every $D > 0$, there exists a partition:
$$
X_k = \bigcup_{1\leq i \leq i_D} X_{k,i}
$$
such that every pair of points $p$, $q$ in the same subset $X_{k,i}$ satisfies $d_k (p,q) \geq D$
and the number $i_D$ of subsets in the partition is a polynomial function of $D$ independent of $k$.

If $p$ lies in $X_{k,i}$, we consider the following partition:
$$
X_k = \{ p \} \cup S_1 \cup S_2 \cup S_3 
$$
where:
\begin{eqnarray*}
S_1 & = & \bigcup_{j<i} X_{k,j}
\\
S_2 & = &  X_{k,i} \backslash \{ p \}
\\
S_3 & = & \bigcup_{j>i} X_{k,j}.
\end{eqnarray*}

Hence:
$$
s_k = \sigma_k ( H_{k,p} , p ) + s_k^1 + s_k^2 + s_k^3
$$
where $s_k^{l} = \sum_{q \in S_l} \sigma_k ( H_{k,q}, q )$.

\medskip 

Let $\varepsilon_1 > \varepsilon_2 > \dots > \varepsilon_{i_D}$ be a sequence of sufficiently small positive numbers.
By Lemma \ref{local avoidance}, for every sufficiently large $k$, there exists a family $(H_{k,p})_{p \in X_k}$ of 
polynomial maps in $\Pol _l (\C^n,\C^m)$ satisfying the following two conditions,
for every $p\in X_{k,i}$:
\begin{enumerate} \item $ \| H_{k,p} \| \leq \varepsilon_i$. 
\item For every $z \in K$:
$$
\dist (A_z , T_l (R_p ( s_k^1 + \sigma_k (H_{k,p}, p)),z)) \geq \varepsilon_i ( - \log \varepsilon_i)^{-N}
$$
where $A_z = \{ H \in \Pol _l (\C^n,\C^m), \; (z,H)\in A \}$ and, as usual, the polynomial map $T_l (R_p ( s_k^1 + \sigma_k (H_{k,p}, p)),z) \in \Pol _l (\C^n,\C^m)$
is the $l-$jet of  the map $R_p ( s_k^1 + \sigma_k (H_{k,p}, p))$ at $z$.

\end{enumerate}

By Lemma \ref{somme de pics}:
$$
\|  T_l ( R_p s_k^3 , z) \| \leq C \varepsilon_{i+1}
$$

Similarly, by the remark following Lemma \ref{somme de pics}:
$$
\|  T_l ( R_p s_k^2 , z) \| \leq C \varepsilon_{i} \exp \left( -  \frac{D^2}{C} \right)
$$
because, for every $q \in X_{k,i} \backslash \{ p\}$:
\begin{eqnarray*}
d_k \left( \varphi_p \left( \frac{z}{ \sqrt{k}}  \right), q \right)
& \geq & d_k (p,q) - d_k \left( \varphi_p \left( \frac{z}{ \sqrt{k}}  \right), p \right)
\\
& \geq & D - C
\end{eqnarray*}
where the constant $C$ is independent of $D$.

\medskip

Since $i_D$ is a polynomial function of $D$,
a calculation shows that for every $C$ and $N$,
if $D$ is sufficiently large then there exists a sequence $\varepsilon_1 > \varepsilon_2 > \dots > \varepsilon_{i_D}$
satisfying the following two estimates:
\begin{eqnarray*}
C \varepsilon_{i+1} & \leq & \frac{1}{4} \varepsilon_i ( - \log \varepsilon_i)^{-N}
\\
C \exp \left( -  \frac{D^2}{C} \right) & \leq & \frac{1}{4} ( - \log \varepsilon_i)^{-N}.
\end{eqnarray*}

Hence, if $p \in X_{k,i}$:

\begin{eqnarray*}
\dist  (A_z , T_l (R_p s_k,z)) 
& \geq &
\dist (A_z , T_l (R_p ( s_k^1 + \sigma_k (H_{k,p}, p)),z))
\\
& - & \| T_l ( R_p s_k^2 , z) \| - \|  T_l ( R_p s_k^3 , z) \|
\\
& \geq & 
\frac{1}{2} \varepsilon_i ( - \log \varepsilon_i)^{-N}.
\end{eqnarray*}

The {\it infimum} of the right-hand side for $1 \leq i \leq i_D$ is positive and independent of $k \gg 1$, $p\in X_k$ and $z \in K$.

\qed

\bigskip

Now, we will complete the proof of Theorem \ref{avoidance sections}.
The trivial holomorphic line budle ${\cal L}: \C \times \C^n \to \C^n$
is endowed with the non-trivial Hermitian metric $\exp \left( - \frac{\pi}{2} \| z \| ^2 \right)$.
Let $u: \C^n \to \C^n$ be an affine unitary transformation:
$$
u(z) = \alpha (z) + \beta.
$$

A lift of $u$ is an automorphism $\widetilde{u}$ of the Hermitan bundle ${\cal L}$:
$$
\widetilde{u} (w,z) = \left(  \lambda w \exp \left( \frac{\pi}{2} \| \alpha ^{-1}( \beta) \|^2 + \pi \langle \alpha ^{-1}( \beta) , z\rangle \right) , u(z)  \right)
$$
where $w\in \C$, $z\in \C^n$
and $\lambda$ is a unit complex number.
Since the Hermitian inner product $\langle .,. \rangle$ is linear in the second argument, the map: 
$$
z \mapsto  \lambda  \exp \left( \frac{\pi}{2} \| \alpha ^{-1}( \beta) \|^2 + \pi \langle \alpha ^{-1}( \beta) , z\rangle \right)
$$
is holomorphic. 
Hence, if a subset $A \subset  \C^n \times \Pol _l (\C^n,\C^m)$ 
satisfies the invariance hypotheses of Theorem \ref{avoidance sections} then, in particular, $A$
is invariant under the action of $\widetilde{u}$ upon $ \C^n \times \Pol _l (\C^n,\C^m)$. 

Let $A$ be a subset of  $\C^n \times \Pol _l (\C^n,\C^m)$ satisfying the hypotheses of Theorem \ref{avoidance sections}
and let $\overline{\B} \subset \C^n$ be the closed unit ball.
We apply Lemma \ref{globalisation} to $A$ and $K = \overline{\B}$.
There exists a renormalizable sequence of holomorphic sections $(s_k)_{k \gg 1}$
such that 
every limit section $s_{\infty}$ of the sequence of renormalized maps $(R(s_k, \tau_{p_k}))_{k\gg 1}$
avoids $A \cap (\overline{\B} \times  \Pol _l (\C^n,\C^m))$, if $p_k$ lies in $X_k$.
Hence, if $z \in \overline{\B}$ then  $T_l ( s_{\infty}, z) \notin A$.

We will prove that a similar result holds if we replace the subset $A \cap (\overline{\B} \times  \Pol _l (\C^n,\C^m))$ with $A$
and if we consider any other normal sequence of framed charts $(\widehat{\varphi}_k, \widehat{\tau}_k)_{k\gg 1}$  
standard at the origin. 

Let $\widehat{s}_{\infty}$ be a limit section of the normal sequence $( R (s_k, \widehat{\tau}_k ))_{k \gg 1}$
and let $\widehat{z}_0\in \C^n$ be a point.
We set $\widehat{p}_k = \widehat{\varphi}_k \left( \frac{\widehat{z}_0}{\sqrt{k}} \right)$.
Since $X_k$ is a $g_k -$discretization, there exists a point $p_k \in X_k$ satisfying $d_k (p_k, \widehat{p}_k) \leq 1$. 
Let $u_k$ be the renormalized change of coordinates:
$$
u_k (z) = \sqrt{k} \, \varphi_{p_k}^{-1} \left( \widehat{\varphi}_k \left( \frac{z}{\sqrt{k}} \right) \right). 
$$
The local biholomorphisms $(u_k)_{k \gg 1}$ form a normal family and every limit $u_{\infty}$
is an affine unitary transformation $\C^n \to \C^n$. 

Let $\widetilde{u}_k$ denote the transition map between
the trivializations
$$z \mapsto \widehat{\tau}_k \left( \frac{z}{\sqrt{k}} \right) \mbox{ and } z \mapsto \tau_{p_k} \left(  \frac{z}{\sqrt{k}}   \right).$$
Similarly, the maps $(\widetilde{u}_k)_{k \gg 1}$ form a normal family and every limit $\widetilde{u}_{\infty}$
is the lift of a limit $u_{\infty}$ of $(u_k)_{k \gg 1}$.

The map $\widehat{s}_{\infty}$ is the limit of a convergent subsequence $( R (s_k, \widehat{\tau}_k ))_{k \in I}$.
We can assume that the corresponding sequences
$( u_k )_{k \in I}$ and $( \widetilde{u}_k )_{k\in I}$
converge to $u_{\infty}$ and $\widetilde{u}_{\infty}$. 
Moreover, since $d_k (p_k, \widehat{p}_k) \leq 1$ for all $k\in I$, we can assume that the points $(u_k (\widehat{z}_0))_{k \in I}$
converge to some point $z_0 \in \overline{\B} \subset \C^n$.

Then the image of $T_l ( \widehat{s}_{\infty}, \widehat{z}_0)$ under the action of $\widetilde{u}_{\infty}$ equals $T_l ( s_{\infty}, z_0)$
where $s_{\infty}$ is a limit map of the sequence $( R (s_k, \tau_{p_k} ))_{k \gg 1}$.
By definition of the sequence $(s_k)_{k \gg 1}$, the jet
$T_l ( s_{\infty}, z_0) $ does not lie in  $A$.
Since $A$ is invariant under the action  of $\widetilde{u}_{\infty}$,
we conclude that $T_l ( \widehat{s}_{\infty}, \widehat{z}_0)  \notin A$.
The sequence $(s_k)_{k \gg 1}$ avoids $A$ asymptotically and the proof of Theorem \ref{avoidance sections} is completed.

\subsection{Proof of the avoidance theorem for submanifolds}\label{proof sub}
Here, we complete the proof of Theorem \ref{avoidance}.

\begin{defn}
For every sufficiently large $k$, let $s_k$ be a section of $\C^m \otimes L^{\otimes k}$
such that the sections $(s_k)_{k\gg 1}$ form a renormalizable sequence.
The sequence $(s_k)_{k\gg 1}$ is {\it transverse to $0$ asymptotically}
if every limit section $s_{\infty}: \C^n \to \C^m$ is transverse to $0$.

\end{defn}

In order to construct submanifolds of dimension $d$, we set $m = n-d$.

\begin{lem}\label{sections submanifolds}
For every sufficiently large $k$, let $s_k$ be a holomorphic section of $\C^{n - d} \otimes L^{\otimes k}$
such that the sections $(s_k)_{k\gg 1}$ form a renormalizable sequence
transverse to $0$ asymptotically.

Then, for every sufficiently large $k$, the section $s_k$ is transverse to $0$ 
and the zero sets $(Y_k = \{ s_k = 0\} )_{k\gg 1}$ form a renormalizable sequence of complete intersections
of dimension $d$.
Moreover, for every limit submanifold $Y_{\infty}$ of $(Y_k )_{k\gg 1}$,
there exists a limit section $s_{\infty}$ of $(s_k)_{k\gg 1}$ such that $Y_{\infty}$ equals $\{ s_{\infty} = 0 \}$,
up to a linear transformation of $\C^n$.

\end{lem}

\proof
First, we will prove by contradiction that, for every sufficiently large $k$, the section $s_k$ is transverse to $0$.
Assume that, for infinitely many $k$, there exists a point $p_k \in \{ s_k = 0 \} $ such that $s_k$ is not transverse to $0$ at $p_k$.
For every sufficiently large $k$, let
$(\varphi_k , \tau_k)$ be a framed chart centered at $p_k$ standard at the origin. 
Assume that the framed charts 
$(\varphi_k, \tau_k)_{k\gg 1}$ form a normal family.
We denote by $Rs_k$ the corresponding renormalized map.
A subsequence $(Rs_k)_{k\in I}$ tends to a map $s_{\infty}$ transverse to $0$.
Hence, for every sufficiently large $k\in I$, the map $Rs_k$ is transverse to $0$, say at the origin.
On the other hand, $s_k$ is not transverse to $0$ at $p_k = \varphi_k (0)$. This is a contradiction.

Now, since $s_k$ is tranverse to $0$, the complete intersection $Y_k = \{ s_k = 0 \}$
is smooth.
We will prove that the submanifolds
$(Y_k)_{k \gg 1}$
form a renormalizable sequence.

For every large $k$, let $\widehat{\varphi}_k:  \B \to X$ be a holomorphic chart such that the charts 
$(\widehat{\varphi}_k)_{k\gg 1}$ form a normal family and satisfy the non-degeneracy condition of Definition \ref{renorm sequence}.

We didn't assume that $d\widehat{\varphi}_k (0)$ is unitary. 
Nonetheless, we write $\widehat{\varphi}_k = \varphi_k \circ u_k$ were $u_k$ is a linear transformation $\C^n \to \C^n$
and $\varphi_k$ is a holomorphic map such that $d \varphi_k (0) $ is unitary.
The non-degeneracy condition implies that we can assume that the sequence $(u_k)_{k\gg 1}$ is relatively compact in $GL (n, \C )$.
Hence,
the chart $\varphi_k$  is defined on an ellipsoid
containing a ball $\B (r)$ whose radius $r$ is independent of $k$
and the charts 
$(\varphi_k : \B (r) \to X)_{k\gg 1}$ form a normal family.
Notice that every limit submanifold $\widehat{Y}_{\infty}$ defined by using the charts $(\widehat{\varphi}_k)_{k\gg 1}$ is the image
 $u_{\infty} (Y_{\infty})$ of a limit submanifold $Y_{\infty}$ defined by using the charts $(\varphi_k)_{k\gg 1}$, where $u_{\infty}$
is a linear transformation $\C^n \to \C^n$.

Now we consider the charts $(\varphi_k)_{k\gg 1}$
and we will show that the corresponding renormalized submanifolds $(RY_k)_{k \gg 1}$ form a normal sequence.

For every $k$, there exist a holomorphic trivialization $(\tau_k: \C \times \B (r) \to \varphi_k^* L)_{k\gg 1}$
such that the framed chart $(\varphi_k , \tau_k)$ is standard at the origin and 
the family $(\varphi_k , \tau_k)_{k \gg 1}$ is normal.
As usual, we denote by $Rs_k: \B (\sqrt{k} r) \to \C^m$ the corresponding renormalized map.

Since the sequence $(s_k)_{k \gg 1}$
is renormalizable, every subsequence of $(Rs_k)_{k \gg 1}$
admits a subsubsequence which converges to a limit map $s_{\infty}: \C^n \to \C^m$.
By assumption, $s_{\infty}$ is transverse to $0$.
Hence the submanifolds $RY_k = \{ Rs_k = 0 \}$ converge to the submanifold $\{ s_{\infty} = 0 \} \subset \C^n$ 
in the smooth compact-open topology.

\qed


\begin{lem}\label{lem A B}
Let $A \subset  \Jet^l_{d,n} $ be a set of $l-$jets of complex submanifolds,
which is a closed complex algebraic subset of $\Jet^l_{d,n}$.
We set $m = n -d$.
Then there exists a set
$B$ of $l-$jets of holomorphic maps $\C^n \to \C^m$,
which is a closed complex algebraic subset of
$\C^n \times \Pol _l (\C^n,\C^m)$ of codimension $\geq \min \{ n+1, n - d + \codim (A) \}$
and satisfies the following property:

\medskip

If $s: (\C^n, p) \to \C^m$ is a germ of holomorphic map whose $l-$jet doesn't lie in $B$
then $s$ is transverse to $0$ and, setting $Y = \{ s = 0 \}$,
the $l-$jet of the submanifold $Y$
at $p$ doesn't lie in $A$.

\bigskip

Moreover, we can assume that:
\begin{enumerate}

\item

If $A$ is invariant under the natural action of affine transformations of $\C^n$ upon $\Jet^l_{d,n}$
then $B$ is invariant under their natural action upon $ \C^n \times \Pol _l (\C^n,\C^m)$.

\item

The subset $B$ is invariant under the action (by multiplication) of the germs of holomorphic functions.

\end{enumerate}

\end{lem}

\proof
Let $E \subset   \C^n \times \Pol _l (\C^n,\C^m)$ be the set of $l-$jets $(z,H)$ whose values $H(z)$ equal $0$.
Hence $E$ is a linear subspace of codimension $m$.
Let $B_0 \subset E$ be the closed algebraic subset of jets such that the tangent map $d H (z)$ is not surjective.
The codimension of $B_0$ in $E$ is $\geq n - m + 1$.

\medskip

An element $[ s ] \in E \backslash B_0$ is the jet of a regular equation and hence, $[ s ]$ defines
a jet of submanifold $\cY ( [ s ] ) \in \Jet^l_{d,n}$.
This map $\cY: E \backslash B_0 \to \Jet^l_{d,n}$
is a regular algebraic submersion.
Therefore, $\cY^{-1}(A) \subset E \backslash B_0$ is a closed algebraic subset
whose codimension in $E \backslash B_0$
equals $\codim (A)$.

We set $B = B_0 \cup \cY^{-1} (A)$. Then $B$ is a closed algebraic subset whose
codimension in $E$ is
$\geq \min \{ n - m +1, \codim (A) \}$.
Hence $B$ has the required codimension in $\C^n \times \Pol _l (\C^n,\C^m)$.

Clearly, the subset $B$ satisfies the second invariance condition and if $A$ is invariant then $B$ also satisfies the first one.

\qed

\bigskip

Now, we are able to complete the proof of Theorem \ref{avoidance}.
Let $A \subset \Jet^l_{d,n}$ be a closed algebraic subset satisfying the hypotheses
of Theorem \ref{avoidance}. In particular, $\codim (A) > d$.

By Lemma \ref{lem A B}, the corresponding subset $B \subset  \C^n \times \Pol _l (\C^n,\C^m)$ satisfies the following condition: 
$$\codim (B) \geq \min \{ n + 1, n - d + \codim (A) \}  =  n + 1$$
as well as the other hypotheses of Theorem \ref{avoidance sections}.
Hence, there exists a renormalizable sequence of sections $(s_k)_{k \gg 1}$ which avoids $B$ asymptotically.
By definition of $B$, the sequence $(s_k)_{k \gg 1}$ is transverse to $0$ asymptotically
and, for every limit section $s_{\infty}$, the jets of the submanifold $\{ s_{\infty} = 0 \}$ do no lie in $A$.
By Lemma \ref{sections submanifolds}, if we set $Y_k = \{ s_k = 0 \}$
then the complete intersections $(Y_k)_{k \gg 1}$
form  a renormalizable sequence satisfying the required conditions
and the proof of Theorem \ref{avoidance} is completed. 


\bigskip

{\it E-mail address:} jean-paul.mohsen@univ-amu.fr

\end{document}